\documentclass[reqno,11pt,a4paper,final]{amsart}
\usepackage[margin=1in]{geometry}
\usepackage{amsmath}
\usepackage{amssymb}
\usepackage{amsthm}
\usepackage{amscd}
\usepackage{mathtools}
\usepackage[T1]{fontenc}
\usepackage[utf8]{inputenc}
\usepackage{bbm}
\usepackage[overload]{textcase}
\usepackage[dvipsnames]{xcolor}
\usepackage{esint}
\usepackage{MnSymbol,wasysym}
\usepackage{dutchcal}
\usepackage{xcolor}

\newcommand{\R}{\mathbb{R}}

\newcommand{\T}{\mathbb{T}}

\newcommand{\norm}[1]{\|#1\|}

\newcommand{\dt}{ dt}

\theoremstyle{definition}

\newtheorem*{kom}{Comment}
\newtheorem{tw}{Theorem}

\newtheorem{defi}{Definition}
\newtheorem*{dow}{Proof}

\newtheorem{lem}{Lemma}
\newtheorem{prop}{Proposition}
\newtheorem{przyk}{Example}

\newtheorem*{dzieki}{Acknowledgement}
\newtheorem*{dow2}{Proof of Theorem \ref{main}}
\newtheorem*{dow1}{Proof of Theorem \ref{existence}}
\newtheorem*{tw*}{Theorem 1}
\newtheorem*{twt}{Theorem 2}

\DeclareMathOperator{\esssup}{ess\,sup}

\DeclareMathOperator{\diverg}{div}

\DeclareMathOperator{\supp}{supp}

\DeclareMathOperator{\tr}{tr}
\DeclareMathOperator{\intel}{int}

\renewcommand{\epsilon}{\varepsilon}
\renewcommand{\phi}{\varphi}

\begin{document}

\title{P\MakeLowercase{arabolic} PDE\MakeLowercase{s on low-dimensional structures}}
\author[Ł\MakeLowercase{ukasz} C\MakeLowercase{homienia}]{Ł\MakeLowercase{ukasz} C\MakeLowercase{homienia}, 18.08.2023\\
	I\MakeLowercase{nstitute of} A\MakeLowercase{pplied} M\MakeLowercase{athematics and} M\MakeLowercase{echanics},\\ 
	U\MakeLowercase{niversity of} W\MakeLowercase{arsaw,} P\MakeLowercase{oland}\\
 \MakeLowercase{Address: }B\MakeLowercase{anacha} 2, 02-097 W\MakeLowercase{arsaw} \\
 \MakeLowercase{e-mail: lukasz.chomienia@uw.edu.pl} 
	}
\maketitle

\begin{abstract}
The paper concerns the theory of parabolic equations on a broad class of closed subsets of Euclidean space possessing a kind of tangent structure. A necessary framework for considering evolutionary problems is developed, and fundamental results about the existence and uniqueness of solutions are established. The presented approach is based on applying the semigroup theory of linear operators combined with special geometric constructions related to the shape of the examined structure. The proposed setting is consistent with the theory of elliptic equations on lower dimensional structures and with the second-order calculus of variations developed by Bouchitt{\'e} and Fragala in \cite{Bouchitte}.  
\end{abstract}



\section{Introduction} 
Interest in analysing low-dimensional structures is driven by the need to model physical phenomena appearing on objects with a non-standard geometry. Such kind of issues arise in various engineering sciences, for example, in the theory of strength of materials \cite{bolb}, in the elasticity theory, in the free material design problem \cite{lew}, or in the modelling of heat flow in thin conductors \cite{Rybka}. 

Several known techniques can be applied to formalise variational problems on irregular domains. For instance, in \cite{fat}, the authors use a method based on considering a sequence of auxiliary variational issues posed on tubular neighbourhoods approximating a given lower dimensional submanifold. The most crucial benefit of using this approach is the possibility of avoiding work with sets of the zero Lebesgue measure. Information about the initial problem is retrieved later by passing with radii of tubular neighbourhoods to zero. Unfortunately, the method demands quite restrictive regularity conditions from the considered low-dimensional structure. Such restrictions rule out applications to many structures that naturally appear in various real-life problems. Later it was observed that the modern measure-theoretic methods could be used to avoid obstacles related to the low regularity of examined sets.
An idea standing behind this approach is to represent a given geometric structure as a Radon measure equipped with a tangent bundle. The measure-theoretic framework used in this paper was originally introduced in papers \cite{first}, \cite{omar} to study problems of minimization of functionals on subsets that are singular with respect to the Lebesgue measure.

The calculus of variations in the low-dimensional setting was continuously developed, giving rise to numerous new concepts and generalizations. In \cite{Bouchitte}, a consistent abstract framework of second-order variational problems was developed. A new notion of a weak elliptic partial differential equation in the low-dimensional setting was recently established in paper \cite{Rybka}. Following this new direction, in our forthcoming paper, we study regularity aspects of weak solutions to elliptic PDEs, particularly how such solutions are related to the domain of the second-order operator introduced in \cite{Bouchitte}.

A novelty of the present paper is in developing the theory of parabolic partial differential equations on irregular low-dimensional structures consistent with the variational and elliptic setting studied previously. New results are obtained both for various strong-form second-order issues and weak formulations. We prove key properties of second-order differential operators and establish the main theorem about the existence and uniqueness of solutions to second-order problems with Neumann initial data. Besides that, we establish the connection between solutions to low-dimensional parabolic and elliptic equations. We also propose some new characterisations of Sobolev-type spaces on low-dimensional structures.  

Let $S = \bigcup_{i=1}^mS_i \subset\R^3,$ where $S_i$ are compact, smooth submanifolds (with boundary) of $\R^3,$ which are pairwise transversal and $\dim S_i \in \{1,2\}.$
The goal of this work is to establish the right meaning of the following formal Neumann boundary problem on the geometric structure $S,$   
\begin{equation}\label{formal}
\begin{matrix}
u_t - \diverg(B\nabla u) = 0 & \text{ in } & S \times [0,T], \\ \\
B\nabla u \cdot \nu = 0 & \text{ on } & \partial S \times [0,T], \\ \\
u = g & \text{ on } & S \times \{0\},
\end{matrix}
\end{equation}
where the matrix of coefficients $B$ satisfies a suitable ellipticity condition, $(B\nabla u \cdot \nu)\lfloor_{\partial S}$ is an appropriate normal derivative and $g$ is a given initial data. To translate the problem into a suitable framework, we associate with the structure $S$ a corresponding singular measure $\mu,$ which encodes the geometry of the set $S.$ By an application of the second-order setting of \cite{Bouchitte}, we define a suitable relaxation of the operator $\partial_t  -L,$ where $Lu=\diverg(B\nabla u)$ is related with the underlying measure $\mu.$ 

The main result of the paper deals with the existence of solutions to the strong-form second-order equation:
\begin{tw*}(Existence and uniqueness of solutions)\\
	Assume that $L_{\mu} = \Delta_{\mu}$ (Definition \ref{lap}). For any given initial data $g \in A_T(\Delta_{\mu})_N$ (see, Definition \ref{veryreg} in Section 4.1) of the low-dimensional parabolic problem (equation \eqref{parabolic}) exists a unique solution $u:[0,T]\to L^2_{\mu}$ in the sense of Definition \ref{parabolic2}. 
\end{tw*}
For definitions of the involved operators (Definition \ref{secondop} and Definition \ref{lap}) we refer to Section 2.3. Equation \eqref{parabolic} and Definition \ref{parabolic2} also can be found in Section 2.3 of the paper.
The proof applies the general semigroup theory of differential operators (for example, see \cite{paz}, Section 1.3) and a specific variant of the Hille-Yosida Theorem about operators generating contraction semigroups. A core of this approach is in avoiding using a resolvent operator, which turns out to be not well-defined due to the lack of invertibility of the considered operator. 
A crucial technical component of the proof is verifying the closedness property of the differential operator. This is expressed by the thesis of the following important 
\begin{twt}\label{main}
	Let $\mu \in \mathcal{S}$ be a low-dimensional structure. The operator $L_{\mu}:D(L_{\mu})_N \to L^2_{\mu}$ is closed in the $L^2_{\mu}$-convergence.
\end{twt}
The precise definition of the considered here class of measures is given in Section 2.1, as well as the definitions of the differential operator $L_{\mu}$ and of its domain are located in Section 2.3 (Definition \ref{secondop} and Definition \ref{0domain}, respectively).  
The second-order functional space in which we conduct our reasoning consists of couples $(u,b),$ where $u$ is a function belonging to some proper subspace of the $\mu$-related Sobolev-type space $H^1_{\mu}$ and $b$ is a certain vector field normal to the structure $S.$ Such normal field is called the Cosserat vector field, and its role is to mimic a normal component of the classical gradient $\nabla u.$  Since for a given $u$ the vector field $b$ is not uniquely determined and there is not enough information about its behaviour, problems arise with control of convergence of sequences of tuples $(u_n,b_n).$ To deal with this issue, for a given pair $(u_n,b_n)$ we introduce a certain procedure of modifying the vector field $b_n.$ This process, based on geometrical constructions, produces a new normal vector field $\widetilde{b}_n,$ which has properties similar to the original field $b_n.$ A fundamental reason for introducing the modified sequence of vector fields $\{\widetilde{b}_n\}$ is that its convergence can be completely controlled in terms of corresponding functions $u_n.$ 

Further results of this paper concern alternative approaches to parabolic issues on low-dimensional structures. We needed very regular initial data to determine the existence of solutions in the strong-form second-order problem (see, Theorem \ref{existence}). To extend the class of accessible initial data, we consider the weak formulation of the problem. Then we derive some information about the regularity of weak solutions. 
The connection between parabolic and elliptic problems is expressed in the theorem showing that a weak solution to the parabolic equation converges, as time goes to infinity, to a solution of the stationary heat equation.

To present that the low-dimensional solutions might differ in an essential way from solutions of classical problems, we present two meaningful examples showing that even in the simple case of the stationary heat equation, quite unexpected phenomena appear. 

The paper is organised as follows. 

In Section 2 we evoke basic notions of the new framework. Among others, we present here: definitions of spaces of functions and the adequate low-dimensional framework, preliminary properties of introduced operators, or notions of solutions to parabolic problems.

Section 3 contains elementary and basic results describing the properties of newly introduced objects. This section includes new and vague characterisations of first and second-order spaces of functions.

Section 4 is devoted to the main theorems of the paper -- the proof of the existence of a semigroup generated by the operator $L$ and the result showing the closedness of the operator $L.$ At the end of this section we present an example of the application of the obtained theorems.      

Section 5 is devoted to weak variants of the parabolic heat equation and the regularity of their solutions. We prove that when passing $t \to \infty,$ parabolic solutions converge to a solution of the stationary equation. This section also contains examples of low-dimensional issues.

\begin{dzieki}
Special thanks go to my supervisor -- Professor Anna Zatorska-Goldstein, for her continuous support and help and to Professor Piotr Rybka for suggesting many interesting ideas. \\
The author was in part supported by the National Science Centre, Poland, by the Grant: 2019/33/B/ST1/00535. 
\end{dzieki}

\section{Preliminaries}
For the convenience of readers, we evoke here the most important concepts. 
Definitions related to the first-order functional setting are taken from \cite{first}, \cite{omar}, \cite{man}. Most of the second-order framework definitions come from \cite{Bouchitte}. Elements of the theory of weak elliptic equations are in a form presented in \cite{Rybka}.

Throughout the rest of the paper, if not specified more precisely, the letter $\mu$ denotes a positive real-valued Radon measure on $\R^3.$

\subsection{Lower dimensional structures}
We introduce a main object of interest -- the class of low-dimensional structures. This class consists of certain Radon measures that represent geometrical structures.

The presented class of measures is a modification of that introduced in papers \cite{first}, \cite{omar}.

Let $\Omega \subset \R^3$ be a non-empty, open, bounded and connected subset. For a fixed $m \in \mathbb{N}$ and fixed $k \in \{1,...,m\},$ for any $i\leqslant k,$ let a closed set $S_i \subset \Omega$ be a 2-dimensional smooth submanifold with boundary embedded in $\R^3,$ and for any $i> k,$ let the closed set $S_i \subset \Omega$ be a 1-dimensional smooth submanifold with boundary embedded in $\R^3.$ Assume that for each $i,j \in {1,...,m},$ $i\neq j$
\begin{itemize}
		\item[a)] $\partial \Omega \cap S_i = \partial S_i,$
		\item[b)] $S_i$ is transversal to $S_j,$
		\item[c)] $S_i\cap S_j \cap S_l = \emptyset,$ for $l\in \{1,...,m\} \setminus \{i,j\}.$
\end{itemize}   
With each $S_i$ we associate the measure $\mathcal{H}^{\dim S_i}\lfloor_{S_i},$ where $\mathcal{H}^{\dim S_i}\lfloor_{S_i}$ is the $\dim S_i$-dimensional Hausdorff measure restricted to $S_i.$

A positive Radon measure $\mu$ belongs to the class $\widetilde{\mathcal{S}}$ if it is of the form $$\mu = \sum_{i=1}^m \mathcal{H}^{\dim S_i}\lfloor_{S_i}.$$

Elements of the set $\widetilde{\mathcal{S}}$ are called low-dimensional structures.
By $\mathcal{S}$ we denote the subclass of the class $\widetilde{\mathcal{S}}$ consisting of low-dimensional structures whose component manifolds have fixed dimension, that is, if $\mu \in \mathcal{S}$ and $\supp\mu = \bigcup_{i=1}^mS_i,$ then for all $i\in \{1,...,m\}$ $\dim S_i = k,$ for $k=1$ or $k=2.$
The subclass $\mathcal{S}$ will be especially important for us when considering problems related to the strong-form operators (see, Section 2.3).

The introduced above class $\widetilde{\mathcal{S}}$ of lower dimensional structures considered in this paper is a proper subset of the one used, for example, in \cite{Rybka}. The difference is in two aspects. First, we assume that a boundary of each component manifold $S_i$ is exactly the set of intersection of $S_i$ with the boundary of the ambient set $\Omega$ (condition a)). In \cite{Rybka}, the authors allow a situation, where $\partial S_i \subset S_j,\; i\neq j.$ A second difference is in an additional restriction on possible types of intersections of components. This is expressed in point c). This restriction rules out common intersections of three or more component manifolds. The above-mentioned conditions are not crucial from our perspective and were introduced to simplify further reasoning and shorten the exposition. We think that methods presented in this paper can be easily modified to cover the wider class of structures such as that considered in \cite{Rybka}.

Let $\mu \in \widetilde{\mathcal{S}}$ and $S = \supp \mu.$ Let $S = \bigcup_{i=1}^m S_i,$ where $S_i$ are component manifolds. We denote 
$$\partial S := \bigcup_{i=1}^m \partial S_i.$$
For any $i \in \{1,...,m\}$ let $n\lfloor_{\partial S_i}:\partial S_i \to \R^3$ denote the uniquely determined outward normal unit vectorfield to $\partial S_i$ in the sense of Spivak \cite{spiv}. 
By the outward normal unit vectorfield to the low-dimensional structure $\mu$ we call the function $n: \partial S \to \R^3$ defined as
$$n:= \bigcup_{i=1}^m n\lfloor_{\partial S_i}.$$

\subsection{First-order space of functions}
To define a derivative of a function given on some (possibly irregular) low-dimensional structure (or in general on a support of a Radon measure), we need to develop a notion of a tangent bundle of a measure. Before formalising this object, let us recall the well-known measure-related Lebesgue space.

	 For $p \in [1,\infty]$ the Lebesgue space $L^p_{\mu}$ is defined as a subspace of the space of $\mu$-measurable functions such that 
	 \begin{itemize}
	 \item[a)] for $p\in [1,\infty)$ the norm $$\norm{f}_{L^p_{\mu}}:=\left(\int_{\Omega} |f|^p d\mu\right)^{\frac{1}{p}},$$
	 \item[b)] for $p=\infty$ the norm
	 $$\norm{f}_{L^{\infty}_{\mu}}:=\esssup_{\mu} |f|,$$ 
	\end{itemize}
	is finite.

A tangent space to a Radon measure $\mu$ can be defined in many non-equivalent ways. An extensive discussion of this subject with a comparison of various definitions can be found in \cite{frag}. We follow the method of construction proposed in \cite{first} and in \cite{omar}. The main advantage of this approach is its simplicity and direct analogy with a classical tangent space to a smooth manifold.

	Let us introduce a set $\mathcal{t}^{\perp}_{\mu}$ of all smooth functions vanishing on $\supp \mu,$ i.e.,
	$$\mathcal{t}^{\perp}_{\mu}:=\left\{v \in C^{\infty}_c(\R^3): v=0 \text{ on } \supp \mu \right\}$$ 
	and the set $\mathcal{T}^{\perp}_{\mu}$ of gradients of such functions: 
	$$\mathcal{T}^{\perp}_{\mu}:= \left\{ w \in C^{\infty}_c(\R^3;\R^3): w=\nabla v \text{ on } \supp \mu \text{ for some } v \in \mathcal{t}^{\perp}_{\mu} \right\}.$$
	It is easy to observe that a set-valued mapping $T_{\mu}^{\perp}: \R^3 \to P(\R^3),$ 
	$$T_{\mu}^{\perp} (x):= \left\{w(x)\in \R^3: w \in \mathcal{T}^{\perp}_{\mu}  \right\}$$ assigns to each $x \in \R^3$ certain linear subspace of $\R^3.$  
	The space $T_{\mu}(x)$ tangent to the measure $\mu$ at a point $x \in \R^3$ is defined as 
	$$T_{\mu}(x):= \left\{ w \in \R^3: w\perp T_{\mu}^{\perp} (x)  \right\},$$  
	the symbol $\perp$ denotes the orthogonality in the Euclidean scalar product in $\R^3.$

With a notion of the tangent structure to $\mu$, we might project the classical gradient onto it to obtain its measure-related counterpart. Strictly speaking,
for $\mu$ almost every $x \in \R^3$ let $P_{\mu}(x):\R^3 \to T_{\mu}(x)$ denote the orthogonal projection onto the tangent space $T_{\mu}.$ Then the tangent gradient of a function $u \in C^{\infty}_c(\R^3)$ is defined for a.e. $x \in \R^3$ as 
$$\nabla_{\mu} u(x):= P_{\mu}(x) \nabla u(x).$$

We are ready to introduce the basic Sobolev-like space $H^1_{\mu}.$
The Sobolev space $H^1_{\mu}$ is defined as a completion of the space $C^{\infty}_c(\R^3)$ in the Sobolev norm 
$$\norm{\cdot}_{\mu}:=\left( \norm{\cdot}^2_{L^2_{\mu}} + \norm{\nabla_{\mu} \cdot}^2_{L^2_{\mu}}\right)^{\frac{1}{2}}.$$

Combining a useful characterisation of the space $H^1_{\mu}$ established in \cite{Bouchitte} (see Lemma 2.2) with the fact that multiplying a measure by a bounded and separated from zero density does not change a tangent structure \cite{Rybka}, we derive following observations expressing relations between classical tangent spaces and spaces tangent to a measure $\mu \in \widetilde{\mathcal{S}}.$\\
Assume, that $\mu \in \mathcal{S},$ $\supp \mu = \bigcup_{i=1}^m S_i.$ A classical tangent structure on the component manifold $S_i$ is denoted by $T_{S_i},$ and a classical tangent gradient on the component $S_i$ is denoted by $\nabla_{S_i}.$ Then 
	\begin{itemize}
		\item[a)] $T_{\mu}(x) = \sum_{i=1}^mT_{S_i}(x)$ for $\mu$-a.e. $x \in \R^3,$
		\item[b)] if $u \in H^1_{\mu},$ then for $i \in \left\{1,...,m\right\}$ $u\lfloor_{S_i} \in H^1(S_i),$
		\item[c)] if $u \in H^1_{\mu},$ then for $i \in \left\{1,...,m\right\}$ $(\nabla_{\mu}u)\lfloor_{S_i} = \nabla_{S_i}u,$
		\item[d)] for a fixed $i\in\left\{1,...,m\right\},$ let $\phi \in C^{\infty}(\R^3),\; \phi=0$ on $\bigcup_{j\neq i}S_j,$ $u \in L^2_{\mu}$ and $u\lfloor_{S_i} \in H^1(S_i),$ then $\phi u \in H^1_{\mu}.$  
	\end{itemize}

It is shown in \cite{Rybka}, Section 2.2, that in a case in which there exist $i,j\in \{1,...,m\}$ such that $\dim S_i \neq \dim S_j,$ where $S_i, S_j$ are component manifolds of the low-dimensional structure $\mu \in \widetilde{\mathcal{S}},$ the classical Poincar{\'e} inequality for functions of the class $H^1_{\mu}$ is not valid. Due to this fact, the authors of the mentioned paper introduce a weaker version of the Poincar{\'e} inequality (\cite{Rybka}, Thm. 2.1).
 Consider a partition $I_1,...,I_d$ of the set of indexes $\left\{1,...,m\right\}=I_1\cup...\cup I_d,$ where sets $I_k,\; k\in \{1,...,d\}$ are pairwise disjoint and the following conditions are satisfied:
\begin{itemize}
		\item[a)] for any $k\in \{1,...,d\}$ the characteristic function $\chi_{\left\{\bigcup_{i\in I_k}S_i\right\}}$ of a sum of component manifolds with indexes in $I_k$ belongs to the kernel of the tangent gradient operator, that is 
		\begin{equation*}
		\chi_{\left\{\bigcup_{i\in I_k}S_i\right\}} \in \ker \nabla_{\mu}
		\end{equation*}
		\item[b)] the set $I_k,\;$ $k\in \{1,...,d\}$ is the largest set of indices with the property a): if $\alpha \in I_k$ and $\widetilde{I}_k:= I_k \setminus \{\alpha\},$ then 
		\begin{equation*}
		\chi_{\left\{\bigcup_{i\in \widetilde{I}_k}S_i\right\}} \notin \ker \nabla_{\mu}.
		\end{equation*}
	\end{itemize}  
	With any element $I_k$ of the partition, we associate the projection 
	\begin{equation*}
	P_ku:=\chi_{\left\{\bigcup_{i\in I_k}S_i\right\}} \strokedint_{\bigcup_{i\in I_k}S_i} u d \mu,
	\end{equation*}
	where $\strokedint_Xfd\mu := \frac{1}{\mu(X)}\int_Xfd\mu.$
	Then for $\mu \in \widetilde{\mathcal{S}}$ exists a set of positive constants $\left\{C_1,...,C_d\right\}$ such that for all $k\in\{1,...,d\}$ and any $u \in H^1_{\mu}$ the following version of the Poincar{\'e} inequality is satisfied
	\begin{equation}\label{weakpoincare}
	\sum_{j\in I_k} \int_{\Omega}|u-P_ku|^2d \mu_j \leqslant C_k \sum_{j\in I_k}\int_{\Omega}|\nabla_{\mu}u|^2d\mu_j, 
	\end{equation}
	for $\mu_j:=\mu\lfloor_{S_j}.$

\subsection{Second-order framework}
This part is focused on defining the operator of a second derivative and corresponding spaces of functions. A general idea is close to the one used in the previous section to construct first-order Sobolev spaces, but now a formal realisation is much more complicated and technical. The first problem is that, in general, even in the case of the smooth function $u,$ the $\mu$-tangent gradient does not belong to $H^1_{\mu}$ space. Thus we cannot apply the operator $\nabla_{\mu}$ to it. A second significant thing is related to a choice of a smooth approximating sequence of a function defined on some low-dimensional structure. It turns out that a choice of an approximating sequence affects a value of a second derivative of the limiting function. To solve this problematic issue, we introduce a certain kind of a normal vector field $b$ and examine the differentiability of $\nabla_{\mu}u+b.$

The presented here exposition of the second-order framework follows paper \cite{Bouchitte}. Such notions were originally introduced to study the calculus of variations on low-dimensional structures in a measure-oriented setting.

For each $u \in C^{\infty}_c(\R^3)$ by $\nabla^{\perp}u$ we denote a $\mu$-a.e. normal component of a gradient, that is\\ ${\nabla u = \nabla_{\mu} u + \nabla^{\perp} u.}$ The symbol $\R^{3\times 3}_{\text{sym}}$ stands for the space of $3 \times 3$ symmetric matrices, and $\nabla^2u$ is a matrix of second partial derivatives of a function $u.$ We introduce a set of triples 
$$\mathcal{g}_{\mu}:= \left\{(u,\nabla^{\perp}u,\nabla^2u): u \in C^{\infty}_c(\R^3)\right\} \subset L^2_{\mu} \times L^2_{\mu}\left(\R^3; T_{\mu}^{\perp}\right) \times L^2_{\mu}\left(\R^3;\R^{3\times 3}_{\text{sym}}\right).$$

In an analogy with the method of defining the space $T_{\mu}$ tangent to a measure $\mu,$ we construct the space $M_{\mu},$ which intuitively can be seen as a second-order counterpart of the space $T_{\mu}.$

Let $\overline{\mathcal{g}_{\mu}}$ denote the closure of the set $\mathcal{g}_{\mu}$ in the space $L^2_{\mu} \times L^2_{\mu}(\R^3; T_{\mu}^{\perp}) \times L^2_{\mu}(\R^3;\R^{3\times 3}_{\text{sym}})$ equipped with the natural product norm. We define the space 
$$\mathcal{m}_{\mu}:= \left\{z \in L^2_{\mu}(\R^3;\R^{3\times 3}_{\text{sym}}): (0,0,z)\in \overline{\mathcal{g}_{\mu}}\right\}.$$ 
By Proposition 3.3 (ii) in \cite{Bouchitte}, for $\mu$-almost every $x\in \R^3$ there exists a $\mu$-measurable multifunction $M^{\perp}_{\mu}:\R^3 \to P(\R^3),$ where $P(\R^3)$ is the power set of $\R^3$ such that 
$$\mathcal{m}_{\mu} = \left\{z \in L^2_{\mu}(\R^3;\R^{3\times 3}_{\text{sym}}): z(x) \in M^{\perp}_{\mu}(x) \text{ for } \mu-a.e.\; x\in \R^3 \right\}.$$ 
For $\mu$-a.e. $x \in \R^3,$ the space of $\mu$-tangent matrices $M_{\mu}$ is defined as 
$$M_{\mu} := \left\{z \in L^2_{\mu}(\R^3;\R^{3\times 3}_{\text{sym}}): z(x) \perp  M^{\perp}_{\mu}(x) \text{ for } \mu-a.e.\; x\in \R^3\right\}.$$
Here the symbol $\perp$ stands for orthogonality in the sense of the standard Euclidean scalar product in $\R^{3\times 3}.$

For a detailed discussion of properties of $M^{\perp}_{\mu},$ $M_{\mu}$ and the space $\mathcal{g}_{\mu},$ see Lemma 3.2. and Proposition 3.3. in \cite{Bouchitte}.

For $\mu$-a.e. $x\in \R^3,$ let $Q_{\mu}(x):\R^{3\times 3}_{\text{sym}} \to M_{\mu}(x)$ be the orthogonal projection. Consider the set $$D(A_{\mu}):=\left\{(u,b)\in L^2_{\mu} \times L^2_{\mu}(\R^3;T_{\mu}^{\perp}): \exists z\in L^2_{\mu}(\R^3;\R^{3\times 3}_{\text{sym}}) \text{ such that } (u,b,z)\in \overline{\mathcal{g}_{\mu}}\right\}.$$\label{dommu}
A normal vectorfield $b$ such that $(u,b) \in D(A_{\mu})$ is called the Cosserat vectorfield of $u \in L^2_{\mu}.$
On the domain $D(A_{\mu})$ we introduce the operator $A_{\mu}$ as $$A_{\mu}(u,b) := Q_{\mu}(z).$$ This operator can be understand as a $\mu$-related matrix of second-order derivatives of a function $u.$

In the next step, we propose generalising a classical tangent Hessian of a function defined on a smooth manifold.

	For $\mu$-a.e. $x \in \R^3$ let $P_{\mu}^{\perp}(x):\R^3 \to T_{\mu}^{\perp}(x)$ denotes an orthogonal projection onto the normal space $T_{\mu}^{\perp}(x).$  
	Put 
	$$D(\nabla^2_{\mu}):= \{u \in H^1_{\mu}: \exists b : (u,b)\in D(A_{\mu})\}$$ 
	and 
	$$D(C) := \{v \in L^2_{\mu}(\R^3;\R^3): P_{\mu}(x)v(x) \in (H^1_{\mu})^3,\; P_{\mu}^{\perp}(x)v(x) \in (H^1_{\mu})^3 \text{ for } \mu-a.e. \; x\in \R^3\}.$$ Let $C: D(C) \to L^2_{\mu},$ be an operator defined as 
	$$Cv := P_{\mu}^{\perp}\nabla_{\mu}(P_{\mu}v) +  P_{\mu}\nabla_{\mu}(P_{\mu}^{\perp}v),$$ 
	and let $T_C: \mathbb{R}^3 \to \mathbb{R}^{3\times 3}$ be a tensor field 
	$$T_C(x)v(x) := (Cv)(x).$$ 
	A right generalisation of a tangent Hessian to the class of lower dimensional structures is an operator $\nabla^2_{\mu}: D(\nabla^2_{\mu}) \to (L^2_{\mu})^{3\times 3}$ given by a formula 
	$$\nabla^2_{\mu}u := P_{\mu}A_{\mu}(u,b)P_{\mu} - T_Cb.$$

The operator $C: D(C)\to L^2_{\mu}$ might be extended to a continuous operator on the whole $L^2_{\mu}$ space (this can be done by the Hahn-Banach Theorem). This fact implies, that the expression $T_Cb$ makes sense for any Cosserat vector field $b$ related to $u \in D(\nabla^2_{\mu}).$

To verify that the operator $\nabla_{\mu}^2$ is properly defined on the space $D(\nabla_{\mu}^2),$ it is necessary to check that its value is independent of a choice of a Cosserat vector field $b.$ This fact was proven in \cite{Bouchitte} in Proposition 3.15.

Let us recall an observation (Prop. 3.10 in \cite{Bouchitte}) describing more straightforwardly a structure of the operator $A_{\mu}$ if $\mu$ is an element of the class $\mathcal{S}.$ 
Firstly please notice that in the considered case, the following inclusion is true 
\begin{equation}\label{prop1}
D(A_{\mu}) \subset \left\{(u,b)\in H^1_{\mu} \times L^2_{\mu}(\R^3;T_{\mu}^{\perp}): \nabla_{\mu}u+b \in H^1_{\mu}\right\}.
\end{equation}
The above inclusion is strict for a generic measure being a member of $\mathcal{S}.$
A second worth-to-observe thing is that if $\mu \in \mathcal{S},$ then the operator $A_{\mu}$ takes a rather expected form and can be expressed as
\begin{equation}\label{prop2}
A_{\mu}(u,b) = \nabla_{\mu}(\nabla_{\mu}u+b).
\end{equation}

We introduce an operator that can be interpreted as an implementation of the operator $\diverg(B \nabla_{\mu} u)$ in the setting of the second-order low-dimensional theory. As the formula is given in terms of the previously proposed operators, it is consistent both with the second-order theory of Bouchitt{\'e} and Fragala \cite{Bouchitte} and with the theory of first-order differential equations established in \cite{Rybka}.

\begin{defi}(Second-order differential operator $L_{\mu}$)\label{secondop}\\
	Put $D(L_{\mu}) := D(\nabla^2_{\mu})$ and let $B$ be a $3 \times 3$ smooth, elliptic and symmetric matrix, that is, $B=(b_{ij})_{i,j\in \{1,2,3\}},\; b_{ij} \in C^{\infty}(\mathbb{R}^3),\; b_{ij} \geqslant c>0,\; B=B^{T},$ and there exists a constant $C>0$ such that for every $\xi \in \mathbb{R}^3,\; C|\xi|^2\leqslant \sum_{i,j=1}^3b_{ij}\xi_i\xi_j.$ We introduce the tangent divergence operator $\diverg_{\mu}:(H^1_{\mu})^3 \to L^2_{\mu},$ $$\diverg_{\mu} ((v_1,v_2,v_3)) := \tr\nabla_{\mu}(v_1,v_2,v_3).$$ The second-order differential operator $L_{\mu} : D(L_{\mu}) \to L^2_{\mu}$ is defined as $$L_{\mu}u := \sum_{i,j=1}^3 b_{ij}(\nabla^2_{\mu}u)_{ij} + \sum_{i=1}^3 (\nabla_{\mu}u)_i \diverg_{\mu}(b_{i1},b_{i2},b_{i3}).$$ 
\end{defi}

We are especially interested in the case where the matrix of coefficients $B$ is the identity matrix.

\begin{defi}(Operator $\Delta_{\mu}$)\label{lap}\\
Let the operator $L_{\mu}$ be as defined above and assume that the matrix $B$ is the identity matrix, i.e. $B=Id.$ In such case the operator $L_{\mu}$ will be denoted by $\Delta_{\mu}$ and its domain by $D(\Delta_{\mu}).$ The operator $\Delta_{\mu}$ can be expressed as  $$\Delta_{\mu}u = \tr \nabla^2_{\mu}u.$$
\end{defi}

A slightly modified version of the domain $D(L_{\mu}),$ which includes information about Neumann boundary conditions, will also be needed. The definition we give might be without any problems expressed in case of any regular enough flow through the boundary $\partial S$, but for clarity of further proceedings, we pose it in the special case of zero flow through the boundary.

\begin{defi}($D(L_{\mu})$ with Neumann boundary conditions)\label{0domain}\\
	Assume that a matrix $B$ corresponds to $L_{\mu}$ as in Definition \ref{secondop}. By $D(L_{\mu})_N \subset D(L_{\mu})$ we denote a subspace that consists of all $u\in D(L_{\mu})$ such that for all component manifolds $S_i,\; i\in \{1,...,m\}$ we have 
	$$(B\nabla_{S_i} u)\lfloor_{\partial S_i} \cdot n\lfloor_{S_i} = 0,$$ 
	almost everywhere with respect to the measure $\mathcal{H}^{\dim S_i - 1}$ on $\partial S_i.$ Here $n$ is the outward normal unit vector to $\partial S.$
\end{defi}

Directly from the definition of $D(L_{\mu})$ it follows, that if $u\in D(L_{\mu}),$ then for all $i\in \{1,...,m\}$ $u|_{S_i}\in H^2(S_i).$ This implies, that a trace of $\nabla_{S_i} u$ for $u \in D(L_{\mu})$ is well-defined. Moreover, please note, that $D(L_{\mu})_N$ is a Banach subspace of the space $D(L_{\mu})$ in the inherited norm.

Analogously the domain of the operator $A_{\mu}$ can be equipped with Neumann boundary conditions.

\begin{defi}($D(A_{\mu})$ with Neumann boundary conditions)\\
The subspace $D(A_{\mu})_N \subset D(A_{\mu})$ is defined as 
$$D(A_{\mu})_N:=\left\{(u,b)\in D(A_{\mu}): u\in D(L_{\mu})_N\right\}.$$
\end{defi}

From a perspective of the research conducted in this paper, an important role will be played by semigroups of linear operators whose evolution is governed by the operator $L_{\mu}.$
 
We say that the operator $L_{\mu}$ generates a semigroup $S:\left[0,+\infty\right)\times L^2_{\mu} \to L^2_{\mu}$ if for all $u \in D(L_{\mu})$ we have $$L_{\mu}u = \lim_{t \to 0_+} \frac{S(t)u - u}{t},$$ where the limit is taken in the $L^2_{\mu}$-norm.

Let $\mu \in \mathcal{S}$ be a low-dimensional structure. Denote $E := \supp \mu;$ $E \subset \Omega.$  
Let $E=\bigcup_{i=1}^mE_i,$ where $E_i$ is a component manifold and $n$ be a normal unit vector field on $\partial E$ directed outward. Let $l^1$ denotes the one-dimensional Lebesgue measure and assume that $g \in L^2_{\mu}.$ Moreover, let $B$ be a matrix of coefficients of the operator $L_{\mu}$ as in Definition \ref{secondop}.

A low-dimensional counterpart of a parabolic problem with Neumann boundary conditions is defined as 
\begin{equation}\label{parabolic}
\begin{matrix}
\partial_tu - L_{\mu}u = 0 &\; \mu\times l^1-\text{a.e. in } E \times [0,T],\\ \\

B\nabla_{\mu}u\cdot n = 0 &\; \text{ on } \partial E \times [0,T],\\ \\

u = g &\; \text{ on } E \times \{0\}.
\end{matrix}
\end{equation}

\begin{defi}(Solution to a parabolic problem)\label{parabolic2} \\
	A solution to the low-dimensional parabolic problem is a function $u$ determined by a semigroup $S: \left[0,T\right] \times L^2_{\mu} \to L^2_{\mu}$ generated by the operator $L_{\mu}.$ Precisely, $u:\left[0,T\right] \to L^2_{\mu}$ is a solution of parabolic Neumann problem \eqref{parabolic} if 
	\begin{itemize}
	\item[a)] for all $t \in \left[0,T\right],$ 
	$$u(t) = S(t)g,$$
	where $S: \left[0,T\right] \times L^2_{\mu} \to L^2_{\mu}$ is a semigroup generated by the operator $L_{\mu}$ and the function $g$ is a given initial data, 
\item[b)]the function $u$ satisfies the equation $$\partial_tu - L_{\mu}u = 0,$$  $\mu\times l^1-\text{almost everywhere in } E \times [0,T],$
 \item[c)] $u$ satisfies zero Neumann boundary condition in a sense formulated below.
	\end{itemize}
	For every $t \in \left[0,T\right]$ and every $i\in \{1,...,m\}$ we have $B\nabla_{\mu}u\cdot n = 0$ almost everywhere on $\partial E_i$ with respect to $\mathcal{H}^1\lfloor_{\partial E_i},$ where $B\nabla_{\mu}u\cdot n:D(L_{\mu}) \to L^2_{\mathcal{H}^1\lfloor_{\partial E_i}} ,$ 
	$$(B\nabla_{\mu}u\cdot n)(x) := B(x)(\nabla_{\mu}u)\lfloor_{\partial E_i}(x) \cdot n\lfloor_{\partial E_i}(x).$$ 
	$(\nabla_{\mu}u)\lfloor_{\partial E_i}$ is a trace of a Sobolev function on the boundary of $E_i$ and $n(x)$ is the outer normal vector (of length one) at a point $x \in \partial E.$   
\end{defi}

\subsection{Weak formulation of parabolic problem} 
Throughout this section, we do not limit our attention to structures consisting of component manifolds of fixed dimension one or two, but we consider general low-dimensional structures $\mu \in \widetilde{\mathcal{S}}.$ Let us remind that on such structures, the Poincar{\'e} inequality is satisfied in the generalized sense as stated in equation \eqref{weakpoincare}.

Our considerations will be held in the following spaces of functions.
 \begin{defi}
Let $\mathring{L^2_{\mu}}:=\{u \in L^2_{\mu}:\int_{\Omega}u d\mu=0\},$ $\mathring{H}^1_{\mu}:= \{u \in H^1_{\mu}:\int_{\Omega}u d\mu=0\}.$ We define 
$$\mathcal{H}:= L^2(0,T;{H}^1_{\mu}),\; \mathring{\mathcal{H}}:= L^2(0,T;\mathring{{H}^1_{\mu}})$$ and 
$$\mathcal{T}:= \{v \in \mathcal{H}: v \in W^{1,2}(0,T;L^2_{\mu}), v(T)=0\},$$ with the zero mean counterpart $$\mathring{\mathcal{T}}:= \{v \in \mathring{\mathcal{H}}: v \in W^{1,2}(0,T;L^2_{\mu}), v(T)=0\}.$$
The spaces $\mathcal{T}, \mathring{\mathcal{T}}$ are equipped with the norm
$$\norm{u}_{\mathcal{T}}:= \left(\norm{u}^2_{L^2H^1_{\mu}}+\norm{u(0)}_{L^2_{\mu}}\right)^{\frac{1}{2}}.$$
Moreover, let $$\mathcal{C}^{\infty}_0:=\{w\in C^{\infty}([0,T];C^{\infty}_c(\R^N)): w(T)=0\},$$ and $$\mathring{\mathcal{C}^{\infty}_0}:=\{w\in C^{\infty}([0,T];C^{\infty}_c(\R^N)): w(T)=0, \int_{\Omega} w d\mu=0\}.$$
\end{defi}

In the next definition, we define operators needed for establishing the weak counterpart of the parabolic problem. Definition \ref{weak3} gives the precise meaning of the considered issue.  
\begin{defi}
Let $B=(b_{ij})_{i,j\in \{1,2,3\}},\; b_{ij} \in C^{\infty}(\mathbb{R}^3),\; b_{ij} \geqslant c>0,\; B=B^{T},$ and assume there exists a constant $C>0$ satisfying for every $\xi \in \mathbb{R}^3,\; C|\xi|^2\leqslant \sum_{i,j=1}^3b_{ij}\xi_i\xi_j.$
We introduce a bilinear form $E:\mathring{\mathcal{H}} \times \mathring{\mathcal{T}} \to \mathbb{R}$ by a formula $$E(u,v):= \int_0^T\int_{\Omega}(B\nabla_{\mu}u(t), \nabla_{\mu}v(t)) - u(t)v'(t)d\mu dt,$$ and a functional $F:\mathring{\mathcal{T}} \to \mathbb{R}$ defined as 
$$F(v):=\int_0^T\int_{\Omega}f(t)v(t)d\mu dt + \int_{\Omega}u_0 v(0)d\mu.$$ Here we assume that $f \in L^2(0,T;{L^2_{\mu}}),\; u_0 \in \mathring{L^2_{\mu}}$ and $v':=\frac{d}{dt}v.$ 
\end{defi}	        	

\begin{defi}\label{weak3}
Let $f \in L^2(0,T;{L^2_{\mu}})$ and $u_0 \in \mathring{L^2_{\mu}}.$ By a weak parabolic problem with the zero Neumann boundary condition we name an issue of finding a function $u\in \mathring{\mathcal{H}}$ satisfying
\begin{equation}\label{parabolicweak}
E(u,\phi) = F(\phi)
\end{equation}	
for all $\phi \in \mathring{\mathcal{C}^{\infty}_0}.$
\end{defi}

\section{Basic facts}
This section contains elementary facts and observations, which will be applied in further reasoning. \\

The operator $A_{\mu},$ despite its non-standard form, shares some properties expected from a reasonable differential operator. From the perspective of our research, the next observation is very important.

\begin{lem}\label{comp}
	\begin{itemize}
		\item[a)] The operator $A_{\mu}: D(A_{\mu}) \to (L^2_{\mu})^{3 \times 3}$ is closed in the sense of the $L^2_{\mu}$-convergence, that is, if $$(u_n,b_n) \in D(A_{\mu}),\;\; u_n \xrightarrow{L^2_{\mu}} u \in L^2_{\mu},\;\; b_n \xrightarrow{L^2_{\mu}} b \in L^2_{\mu}(\R^3;T^{\perp}_{\mu})$$ and $$A_{\mu}(u_n,b_n) \xrightarrow{L^2_{\mu}} D \in (L^2_{\mu})^{3 \times 3},$$ then $$(u,b) \in D(A_{\mu})\;  \text{  and  }\;  D = A_{\mu}(u,b).$$
		\item[b)] The same is true if the domain $D(A_{\mu})$ is replaced with $D(A_{\mu})_N.$
	\end{itemize}
	\begin{dow}
		\begin{itemize}
			\item[a)] For the proof, see \cite{Bouchitte} Prop. 3.5 (i).
			\item[b)] By the result of point a), the inclusion $D(A_{\mu})_N \subset D(A_{\mu}),$ by applying the classical Poincar{\'e} inequality on each component manifold $E_i,\; (\supp \mu = \bigcup_{i=1}^mE_i),$ and due to continuity of the trace operator it can be checked that the space $D(A_{\mu})_N$ is closed with respect to the convergence described in point a). Thus the proposed result follows. \qed
		\end{itemize} 
	\end{dow}
\end{lem}

Next lemma shows that a global continuity of a function $u \in D(A_{\mu})$ can be derived not only in the case of component manifolds of a fixed dimension but also if the dimensions of components vary.

\begin{lem}\label{conti}
    Let $u \in D(A_{\mu}).$ Assume that $\mu \in \mathcal{S},\; P:= \supp \mu = E_i \cup E_j,\; E_i\cap E_j\neq \emptyset,\; \dim E_i=1 \text{ and } \dim E_j=2.$
	Then the function $u$ is continuous on the set $P.$   
	\begin{dow}
	 To avoid discussion of technical aspects let us assume that $u\lfloor_{E_j}$ is compactly supported in the intel of the component $E_j.$	 
   By inclusion \eqref{prop1} we know that $u\in H^1_{\mu}$ and $u\lfloor_{E_k} \in H^2(E_k)$ for $k=i,j.$ We justify continuity of $u$ proceeding by a contradiction. Suppose our claim is false, i.e. $u$ is discontinuous. As $u\lfloor_{E_k} \in H^2(E_k)$ for $k=i,j$ provides continuity on each component manifold $E_k,\; k=i,j,$ a discontinuity have to appear in the junction set $E_i \cap E_j.$ Let $p \in E_i \cap E_j$ be a fixed point of discontinuity. After multiplying the function $u$ by a constant, we can assume that a jump at the discontinuity point $p$ is greater than $2.$ Let $\{\phi_n\}_{n\in \mathbb{N}},\; \phi_n \in C^{\infty}_c(\R^3)$ be an approximating sequence of $u,$ which justifies membership in $D(A_{\mu}).$ 
		A convergence in the sense of $D(A_{\mu})$ implies convergence in the $H^1$-norm on the 1-dimensional component $E_i$ (of course $D(A_{\mu})$-convergence implies $H^2(E_i)$-convergence, but the weaker type is enough at this moment of reasoning). This implies that it is possible to extract a subsequence converging uniformly on $E_i$ (relabelling is omitted). 
	Let us consider a sequence $w_n:E_j \to \R,$ $w_n:= u\lfloor_{E_j}-\psi_n.$ The sequence $w_n$ satisfies: $w_n \in C_c(\intel E_j), w_n \in H^2(E_j)$ and $w_n \to 0$ in $H^2(E_j).$ Moreover, $|w_n(p)|>1.$ This means that the sequence $w_n$ can be used as a witness showing that the $H^2$-capacity of the point $p$ in the disc $E_j$ is zero. This is a contradiction with the general theory of Sobolev capacity.	
 A general situation where $u\lfloor_{E_j}$ is not necessarily compactly supported in the intel of $E_j$ can be reduced to the discussed above case by modifying the function $u\lfloor_{E_j}$ properly.
		\qed
	\end{dow}
\end{lem}

We introduce a local characterisation of low-dimensional Sobolev-type spaces in terms of the behaviour of functions in a neighbourhood of junction sets. This will allow us to divide our study into a certain number of ``generic'' local cases.

\begin{defi}(Partition of unity)\\
	Let $\mu$ be an arbitrary measure of the class $\mathcal{S}.$ Denote $S:=\supp \mu \subset \mathbb{R}^3,$ where $S=\bigcup_{i=1}^m S_i$ and $S_i$ are component manifolds. 
	We introduce the set
	$$\widetilde{J}:= \left\{p\in S: \exists i,j\in \{1,...,m\},\; i\neq j,\; p\in S_i\cap S_j\right\}.$$ 
	On the set $\widetilde{J}$ we define a relation $\sim_r$ as
	$$p\sim_r q \iff \exists i,j \in \{1,...,m\},\; i\neq j,\; \exists w\in C^{\infty}\left([0,1];S_i \cap S_j\right):\; w(0)=p,\; w(1)=q.$$
	It is easy to check that $\sim_r$ is an equivalence relation and the number of all equivalence classes is finite. In each equivalence class $[p]_r$ we arbitrarily choose one representative $p.$ Let $J$ stands for the set of all such representatives. Every element of the set $J$ represents a different junction set of the structure $S.$ 
	
        A standard result of the mathematical analysis determines the existence of a finite set $K$ and a partition of unity that consists of functions $\alpha_p \in C^{\infty}_c(\R^3),\; p\in I:=J\cup K$ and two sequences of sets (open in the topology of $\R^3$): 
	$\left\{O_p\right\}_{p \in I},\; \left\{U_p\right\}_{p \in I} \subset \R^3$ satisfying stated below properties.
	For each $p\in J,\; p \in S_i\cap S_j$ exists an open set $O_p \supset S_i\cap S_j$ and exists an open set $U_p$ such that $S_i \cap S_j \subset U_p\subset O_p,$ $U_p \cap O_q = \emptyset,$ for $q\in I,\; q\neq p.$  The family $\left\{O_p\right\}_{p \in I}$ cover the set $S,$ that is $\bigcup_{p \in I} O_p \supset S.$ For each $p\in I$ the functions $\alpha_p$ satisfy $\alpha_p \subset \overline{O_p}$ and $\sum_{p\in I}\alpha_p  \equiv 1$ on $S.$ We also assume that the low-dimensional structure $\mu \in \mathcal{S}$ satisfies two additional facts. For each $p\in J$ there exists $x \in S_i \cap S_j$ and a ball $B(x,1)$ such that $B(x,1) \subset S$ and $O_p = B(x,1).$ Moreover, the family $\{O_p\}_{p\in I}$ can be chosen in a way that for each $O_p$ we have $\partial O_p \in C^{2}.$
\end{defi}

\begin{kom}
The last two assumptions given on the measure $\mu \in \mathcal{S}$ exclude from our considerations cases of some low-dimensional structures. These two additional conditions can be cancelled without difficulties, resulting in extended computations and a more complicated exposition of further proofs.   
\end{kom}

The following propositions show that being a member of the spaces $H^1_{\mu}$ or $D(A_{\mu})$ can be completely characterised in terms of behaviour near points of junctions.

\begin{prop}(``from local to global'' characterisation of $H^1_{\mu}$)\label{loc1}\\
	Denote $\mu_p := \mu\lfloor_{O_p}, p\in I.$ The space $H^1_{\mu}$ can be characterised as 
	$$H^1_{\mu} = \left\{u \in L^2_{\mu}: \forall p\in I\; \forall i\in \{1,...,m\}\; \exists u_p \in H^1_{\mu_p},\; u = \sum_{p \in I}\alpha_pu_p\right\}.$$
	\begin{dow}
		$(\implies)$ Assume that $u \in H^1_{\mu},$ put $u_p:= u\lfloor_{\supp \alpha_p}$ for $p\in I.$
		By standard properties of a partition of unity, it follows that $u = \sum_{p \in I}\alpha_pu_p.$
		
		$(\Longleftarrow)$ Let $u \in \left\{w \in L^2_{\mu}: \forall p\in I\;  \exists w_p \in H^1_{\mu_p}\; w = \sum_{p \in I}\alpha_pw_p\right\}.$ 
		For any $u_p \in H^1_{\mu_p},$ the functions $\alpha_p u_p$ can be extended by zero to obtain $\alpha_p u_p\in H^1_{\mu}.$ Let $\phi^n_p \in C^{\infty}(\R^3)$ be a sequence witnessing belonging of $u_p$ to $H^1_{\mu_p}.$
		Taking the sequence 
		$$\sum_{p \in I}\alpha_p \phi^n_p \in C^{\infty}_c(\R^3)$$ 
		and passing to the limit $n \to \infty$ we conclude that $u \in H^1_{\mu}.$
		\qed
	\end{dow}
\end{prop}

\begin{prop}(``from local to global'' characterisation of $D(A_{\mu})$)\label{propp}\\
	 Denote $\mu_p := \mu\lfloor_{O_p}, p \in I.$  Then
	\begin{equation}\label{loc2}
	D(A_{\mu}) = \left\{(u,b)\in L^2_{\mu}\times L^2_{\mu}(\R^3;T^{\perp}_{\mu}): \forall p\in I\; \exists (u_p,b_p)\in D(A_{\mu_p}),\; u = \sum_{p \in I}\alpha_p u_p\right\}.
	\end{equation}
	\begin{dow}
		$(\implies)$ Let $(u,b) \in D(A_{\mu}).$ For each $p \in I,$ take $u_p:= \alpha_p u.$ By a definition of the domain $D(A_{\mu})$ there exists a sequence $\psi_n \in C^{\infty}(\R^3),$ such that $(\psi_n,\nabla^{\perp}\psi_n) \to (u,b)$ in the sense of the $L^2_{\mu}\times L^2_{\mu}(\R^3;T^{\perp}_{\mu})$-convergence and $A_{\mu}(\psi_n,\nabla^{\perp}\psi_n) \to A_{\mu}(u,b)$ as $n \to \infty.$ For any $p \in I$ consider sequence $\alpha_p\psi_n.$ It is easy to see, that there exist $b_p \in L^2_{\mu_p}(\R^3;T^{\perp}_{\mu_p})$ such that
		$$A_{\mu_p}(\alpha_p \psi_n,\nabla^{\perp}(\alpha_p \psi_n)) \to A_{\mu_p}(\alpha_p u, b_p)$$
		in the norm of $\left(L^2_{\mu}\right)^{3 \times 3},$ when we pass with $n \to \infty.$
		This shows that 
		$$(\alpha_p u, b_p) \in D(A_{\mu_p}) \text{ for all } p\in I$$ 
		and proves the considered implication.
		
		$(\Longleftarrow)$ Assume that $(u,b)$ belongs to the right hand side of equality \eqref{loc2}.
		Let $\psi^n_p \in C^{\infty}_c(\R^3)$ be a sequence justifying membership $(u_p,b_p) \in D(A_{\mu_p}).$ Let take $\alpha_p$ and consider $\alpha_p \psi^n_p$ for $p \in I.$
		The functions $\alpha_p u_p,$ can be extended by zero to the whole set $\supp \mu,$ thus $(\alpha_p u_p,\nabla^{\perp}(\alpha_p u_p))\in D(A_{\mu}).$
		Each term of the sum 
		$$\sum_{p\in I}\alpha_p \psi^n_p \in C^{\infty}(\R^3)$$  
		converges in the sense of $D(A_{\mu}),$ thus $(u,b)\in D(A_{\mu}).$
		\qed
	\end{dow}
\end{prop}

\begin{prop}(``from local to global'' characterisation of $D(A_{\mu})_N$)\label{proppp}\\
	Put $\mu_p := \mu\lfloor_{O_p}, p \in I.$ Then
	
	\begin{equation}\label{loc3}
	D(A_{\mu})_N = \left\{(u,b)\in L^2_{\mu}\times L^2_{\mu}(\R^3;T^{\perp}_{\mu}): \forall p\in I\;  \exists (u_p,b_p)\in D(A_{\mu_p})_N,\; u = \sum_{p \in I}\alpha_p u_p \right\}.
	\end{equation}
	
	\begin{dow}
		The proof follows exactly the same lines as the proof of Proposition \ref{propp}. The needed modification is to change for all $p \in I$ the domains $D(A_{\mu_p})$ to its counterparts equipped with the zero Neumann boundary conditions -- the spaces $D(A_{\mu_p})_N.$
		\qed
	\end{dow}
\end{prop}

\section{Main results}
This section aims to establish the existence and uniqueness of solutions to parabolic equations posed on lower dimensional structures. Precisely we will prove the following theorem:
\begin{tw}\label{existence}(Existence and uniqueness of solutions)\\
	Assume that $L_{\mu}=\Delta_{\mu}.$ For any given initial data $g \in A_T(\Delta_{\mu})_N$ (see, Definition \ref{veryreg} in Section 4.1) of the low-dimensional parabolic problem (equation \ref{parabolic}) exists a unique solution $u:[0,T]\to L^2_{\mu}$ in the sense of Definition \ref{parabolic2}. 
\end{tw}
To prove this result, we will adapt the semigroup theory, and a conclusion will be derived by applying a certain variant of the Hille-Yosida Theorem, treating semigroups generated by bounded linear operators. Before we proceed to the proof of Theorem \ref{existence}, we need to prove Theorem \ref{main}, which is the main technical result of this paper. 
\begin{tw}\label{main}
	Let $\mu \in \mathcal{S}$ be a low-dimensional structure. The operator $L_{\mu}:D(L_{\mu})_N \to L^2_{\mu}$ is closed in the $L^2_{\mu}$-convergence.
\end{tw}
It turns out that the proof of this property in the case of the operator $L_{\mu}$ is more demanding than showing the closedness of $A_{\mu}.$ This phenomenon is related to the fact that a convergence of functions in the sense of $D(L_{\mu})$ space does not provide control of Cosserat vector fields.

\begin{dow2}
Let $\{u_n\}_{n \in \mathbb{N}} \subset D(L_{\mu})_N$ be a sequence such that $u_n \xrightarrow{L^2_{\mu}} u \in L^2_{\mu}$ and $L_{\mu}u_n \xrightarrow{L^2_{\mu}} B \in L^2_{\mu}.$ We need to show that $u \in D(L_{\mu})_N$ and $B = L_{\mu}u.$ A proposed strategy of the proof is based on the fact that the operator $A_{\mu}$ is closed (see Lemma \ref{comp}). For the sequence $\{u_n\}_{n \in \mathbb{N}} \subset D(L_{\mu})_N$ we modify the corresponding Cosserat sequence $\{b_n\}_{n \in \mathbb{N}}$ and construct a new sequence of normal vector fields $\{\widetilde{b}_n\}_{n \in \mathbb{N}}$ convergent in the $L^2_{\mu}$-norm and for which the result of Lemma \ref{comp} is valid.

As $D(L_{\mu})_N \subset D(L_{\mu})$ for any element $u_n \in D(L_{\mu})_N$ exists a sequence $\{u^m_n\}_{m \in \mathbb{N}} \subset C^{\infty}_c(\mathbb{R}^3)$ such that $u^m_n \xrightarrow{L^2_{\mu}} u_n,\; \nabla^{\perp}u^m_n \xrightarrow{L^2_{\mu}} b_n$ and $A(u^m_n, \nabla^{\perp}u^m_n) \xrightarrow{L^2_{\mu}} A(u_n, b_n).$ Existence of such sequences is ensured by the definition of the domain $D(A_{\mu})$ (Definition \ref{dommu}). After extracting, if it is necessary, a subsequence from the sequence $\{u^m_n\}_{m \in \mathbb{N}}$ (we omit to relabel of the chosen subsequence) we may assume that  
\begin{equation}\label{esty2}
\norm{u^m_n-u_n}_{L^2_{\mu}}\leqslant \frac{1}{m^2},\;\;\;\; \norm{\nabla^{\perp} u^m_n-b_n}_{L^2_{\mu}}\leqslant \frac{1}{m^2},\;\;\;\; \norm{A_{\mu}(u^m_n,\nabla^{\perp} u^m_n)-A_{\mu}(u_n,b_n)}_{L^2_{\mu}}\leqslant \frac{1}{m^2}.
\end{equation}
By the weaker Poincar{\'e} inequality \eqref{weakpoincare} or alternatively by using the classical Poincar{\'e} inequality on each component manifold separately we notice that $\norm{\nabla_{\mu} u^m_n - \nabla_{\mu} u_n}_{L^2_{\mu}} \leqslant \frac{1}{m^2}$ possibly after passing to a subsequence once again.

Before we start the process of modifying the Cosserat vectors $b_n,$ we will show that without any loss on generality, the domain can be 'straightened out' locally.\\ 

From this moment until the end of the proof, we assume that the low-dimensional structure $S$ consists of exactly two component manifolds with a non-trivial intersection. At the end of the proof, we will evoke Proposition \ref{proppp} to conclude the global instance. We consider two fixed manifolds $E_i,E_j$ such that $E_i \cap E_j \neq \emptyset$ and treat separately each of the following instances: $\dim E_i=\dim E_j=2,$ $\dim E_i=\dim E_j=1,$ and $\dim E_i=2,\;\dim E_j=1.$\\ 

Let $\mu \in \mathcal{S},$ $\supp \mu = E_i \cup E_j,$ $E_i\cap E_j \neq \emptyset,$ $\dim E_i = \dim E_j = 2.$ Let us assume that each component $E_k, k=i,j$ can be parametrized by a single parametrization -- existence of a partition of unity similar to the one presented in Proposition \ref{proppp} justifies, that such restriction do not narrow the class of considered structures.\\
Let $\Psi_i$ be a parameterization of $E_i$ (up to the boundary) and $\Psi_i(B^1)=E_i,$ where $B^1$ is a 2-dimensional unit ball in variables $(x,y)$ with the center at zero.  We define a diffeomorphism 
$\Theta_i$ by the formula $\Theta_i(x,y,z):= (0,0,z)+\Psi_i(x,z).$ By a construction of $\Theta_i$ we have $\Theta_i^{-1}(E_i \cup E_j) = B^1 \cup F_j,$ where $F_j$ is some 2-dimensional manifold and $B^1 \cap F_j$ is a smooth curve. Let $\Psi_j$ be a parametrisation of $F_j.$ Existence of such mappings is provided by assumptions posed on component manifolds in the definition of the class of low-dimensional structures $\mathcal{S}.$ Analogously to what was made in the case of the component $E_i$ we define a diffeomorphism $\Theta_j$ by the formula $\Theta_j(x,y,z):=(0,y,0) + \Psi_j(x,y).$ Now $\Theta_j^{-1}(B^1\cup F_j) = P^1 \cup P^2,$ where $P^1,P^2$ are ``flat'' 2-dimensional manifolds with $P^1 \cap P^2 \subset \left\{(x,y,z)\in \R^3: y=z=0\right\}.$ Let us denote $I:= \Theta_i\circ \Theta_j,$ thus clearly $I^{-1}(E_i \cup E_j) = P^1 \cup P^2.$ Having the ``flat'' structure $P^1 \cup P^2$ we can easily find a smooth diffeomorphism mapping both components $P^1$ and $P^2$ to $2$-dimensional discs. 
If $\mu \in \mathcal{S}$ is the given low-dimensional structure and $\nu := \mathcal{H}^2\lfloor_{P^1}\otimes \mathcal{H}^2\lfloor_{P^2},$ then obviously $\nu \in \mathcal{S}$ and the $D(A_{\mu})$-convergence is equivalent to the convergence in the sense of $D(A_{\nu}).$ To see this, let us notice that the pushback of the measure $\mu$ by the diffeomorphism $I^{-1}$ is a measure absolutely continuous with respect to $\nu,$ and its $\nu$-related density will be bounded and separated from zero.
A simple computation shows that convergence in the sense of $H^1_{\nu}$ is equivalent to $H^1_{\widetilde{\mu}}$-convergence, $D(A_{\nu})$-convergence is equivalent to $D(A_{\widetilde{\mu}})$-convergence. This also implies that convergence in the sense of $D(L_{\nu})$ is equivalent to convergence in the sense of $D(L_{\widetilde{\mu}}).$ An important consequence of this observation is that in further studies it is enough to focus our attention on the case of the measure $\nu$ being a sum of Hausdorff measures representing each ``flattened'' component manifold.   \\

$\boldsymbol{\dim E_i=\dim E_j =2}$

We start by examining a case of $\dim E_i=\dim E_j=2.$ Further considerations can be conducted in the coordinate system related to components $E_i,\;E_j.$ This means that without loss of generality we can take $$E_i = \left\{(x,y,0)\in \R^3: x^2+y^2\leqslant 1,\; x,y\in [-1,1]\right\},$$ $$E_j = \left\{(x,0,z)\in \R^3: x^2+z^2\leqslant 1,\; x,z\in [-1,1] \right\}$$ and 
$$\mu = \mathcal{H}^2\lfloor_{E_i} \otimes \mathcal{H}^2\lfloor_{E_j}.$$

As the intersection $E_1\cap E_2$ is the set $\left\{(x,0,0)\in \R^3: x\in [-1,1]\right\}$ we know that for any function $w \in D(L_{\mu})$ a first coordinate of the gradient vector $\nabla_{\mu} w = (w_{x,\mu}, w_{y,\mu}, w_{z,\mu})$ belongs to $H^1_{\mu},$ that is $w_{x,\mu}\in H^1_{\mu}.$ The Cosserat vector field $b$ is needed only to ensure that $(0,w_{y,\mu}, w_{z,\mu}) + b \in (H^1_{\mu})^3.$

For $n \in \mathbb{N},$ define 
$$h_{n,\mu}=\left(0,{h_{n1}}_{\mu},{h_{n2}}_{\mu}\right) := \left(0,{u_{n+1}}_{y,\mu},{u_{n+1}}_{z,\mu}\right) - \left(0,{u_n}_{y,\mu},{u_n}_{z,\mu}\right) \in \left(H^1(E_k)\right)^3,\; k=i,j.$$ 
By a fact that $L_{\mu}u_n,\; n \in \mathbb{N}$ is a Cauchy sequence in $L^2_{\mu}$ we are going to show that, after extracting a suitable subsequence, we have for all $n \in \mathbb{N}$
$$\norm{h_{n,\mu}}_{H^1(E_k)}\leqslant \frac{1}{n^2},\; k\in \{i,j\}.$$ 
The operator $L_{\mu}$ can be decomposed into two classical operators acting on each component $E_k,\; k=i,j$ separately. This observation allows us to apply to restricted operators classical estimates of the elliptic regularity theory (with Neumann boundary conditions). For the detailed formulation of the recalled facts, see the appendix in \cite{lasica}, or \cite{Evans} Section 6.3.2 for the Dirichlet boundary condition version. On each $E_k$ this gives the estimate 
\begin{equation}\label{esty}
\norm{h_{n,\mu}}_{H^1(E_k)}\leqslant C\left(\norm{L_{\mu}u_{n+1} - L_{\mu}u_n}_{L^2(E_k)} + \norm{u_{n+1} - u_n}_{L^2(E_k)}\right).
\end{equation}
By the $L^2_{\mu}$-convergence of the right-hand side, after passing to a subsequence if necessary, we obtain that for $k=i,j\;$ $\norm{h_{n,\mu}}_{H^1(E_k)} \leqslant \frac{1}{n^2}.$

Further, for each $n,m \in \mathbb{N}$ the function $h^m_{n,\mu}$ is introduced as $$h^m_{n,\mu}=\left(0,{h^m_{n1}}_{\mu},{h^m_{n2}}_{\mu}\right) := \left(0,{u^{m+1}_{n+1}}_{y,\mu},{u^{m+1}_{n+1}}_{z,\mu}\right) - \left(0,{u^m_n}_{y,\mu},{u^m_n}_{z,\mu}\right) \in \left(H^1(E_k)\right)^3,\; k=i,j.$$ 
Applying the triangle inequality, estimate \eqref{esty} for the $H^1(E_k)$-norm of $h_{n,\mu},$ estimates \eqref{esty2} for the speed of convergence of sequences of smooth functions $A_{\mu}(u^m_n,\nabla^{\perp}u^m_n)$ and $u^m_n,$ combined with the Poincar{\'e} inequality on each component manifold $E_k,$ the following estimate is obtained (for $k \in \{i,j\}$) 
\begin{equation}
\begin{aligned}\label{hestimate}
\norm{h^m_{n,\mu}}_{H^1(E_k)}\leqslant \norm{u^{m+1}_{n+1} - u^m_{n+1}}_{H^2(E_k)} +\norm{u^m_{n+1} - u_{n+1}}_{H^2(E_k)} \\+ \norm{h_{n,\mu}}_{H^1(E_k)} + \norm{u^m_{n} - u_{n}}_{H^2(E_k)} \leqslant \frac{1}{n^2} +\frac{3}{m^2}.
\end{aligned}
\end{equation}

For each smooth function $u^m_n \in C^{\infty}_c(\R^3)$ we know that $(u^m_n, \nabla^{\perp}u^m_n) \in D(A_{\mu}).$ To ensure good behaviour of the Cosserat sequence our next goal is to modify the sequence of pairs $(u^m_n,b^m_n)\in D(A_{\mu}),$ where $b^m_n = \nabla^{\perp}u^m_n.$ Precisely speaking, for each $u^m_n,$ we construct a new normal vector field $\widetilde{b^m_n},$ for which $(u^m_n, \widetilde{b^m_n})$ is also a member of $D(A_{\mu})$ and a sequence of diagonal elements of the presented below infinite lower triangular matrix\\ 
\begin{center}
$\begin{bmatrix}
\nabla_{\mu} u^1_1 + \widetilde{b^1_1} &  &  &  \\
\nabla_{\mu} u^2_1 + \widetilde{b^2_1} & \nabla_{\mu} u^2_2 + \widetilde{b^2_2} &  &  \\
\nabla_{\mu} u^3_1 + \widetilde{b^3_1} & \nabla_{\mu} u^3_2 + \widetilde{b^3_2} & \nabla_{\mu} u^3_3 + \widetilde{b^3_3} &\\
\vdots & \vdots & \vdots & \ddots
\end{bmatrix}.$
\end{center}
converges in the space $(H^1_{\mu})^3.$\\
We are going to show that the diagonal sequence $\left((u^n_n,\widetilde{b^n_n})\right)_{n \in \mathbb{N}}$ is such that $(u^n_n,\widetilde{b^n_n})\in D(A_{\mu}),\;$ $u^n_n \xrightarrow{L^2_{\mu}} u,\;$ $\widetilde{b^n_n} \xrightarrow{L^2_{\mu}} \widetilde{b}$ for some $u, \widetilde{b} \in L^2_{\mu},$ and $\nabla_{\mu}\left(\nabla_{\mu} u^n_n+\widetilde{b^n_n}\right) \xrightarrow{L^2_{\mu}} \nabla_{\mu}\left(\nabla_{\mu} u+\widetilde{b}\right) \in L^2_{\mu}.$

Let a restriction of $w \in L^2_{\mu}$ to the component $E_k,\; k\in \{i,j\}$ be denoted by $w^k,$ that is $w^k:= w\lfloor_{E_k}.$
For the function $u^m_n$ the new Cosserat vectorfield $\widetilde{b^m_n}\in L^2_{\mu}(\R^3;T^{\perp}_{\mu})$ is defined as  
$$\widetilde{b^m_n}(x,y,z) := \begin{cases}
\left(0,\; 0,\; {u^m_{n\; z}}^j(x,y)\right) & \text{ on }\; E_i\\
\left(0,\; {u^m_{n\; y}}^i(x,z),\; 0\right) & \text{ on }\; E_j
\end{cases}
.$$
An idea standing behind the sequence $\widetilde{b^m_n}$ is to 'copy' the function ${u^m_{n\; y}}^i$ from the component manifold $E_i$ (in variables $(x,y)$) to the component manifold $E_j$ (in variables $(x,z)$) and analogously to 'copy' ${u^m_{n\; z}}^j$ from $E_j$ (variables $(x,z)$) to $E_i$ (variables $(x,y)$).

To ensure that the constructed pairs belong to the domain $D(A_{\mu})$ we use the fact that in the case of $\mu \in \mathcal{S}$ there exists a characterisation of membership in $D(A_{\mu})$ stated in Proposition 3.11 in \cite{Bouchitte}. The mentioned proposition is based on the idea of constructing a smooth approximating sequence of the class $C^{\infty}(\R^3)$ by extending a sequence initially defined on $\supp \mu$ to the whole space $\R^3$ by the help of the Whitney Extension Theorem. To conclude that each pair ${(u^m_n,\widetilde{b^m_n}) \in H^1_{\mu} \times L^2_{\mu}(\R^3;T^{\perp}_{\mu})}$ is an element of the domain $D(A_{\mu})$ we need to verify if the following three conditions are satisfied by
\begin{itemize}
\item[a)] For $k \in \{i,j\},$ we need to check if ${u^m_n}^k \in H^2(E_k)$ and ${b^m_n}^k \in \left(H^1(E_k)\right)^3.$ It turns out to be trivial due to the fact that $u^m_n$ is a smooth function in $\R^3$ and due to the way in which we constructed the Cosserat field $\widetilde{b^m_n}.$
\item[b)] We can use the same arguments as before to verify that the junction set $E_i \cap E_j$ belongs to continuity points of $(u^m_n, \nabla_{\mu}u^m_n + b^m_n).$ This is the second assumption, which needs to be satisfied.
\item[c)] The third condition demands that we check if on each component $E_k$ the corresponding restriction ${u^m_n}^k$ is of class $C^2(E_k),$ the normal vector field ${\widetilde{b^m_n}}^k$ is Lipschitz continuous on $E_k$ and $({u^m_n}^k, \nabla_{E_k}{u^m_n}^k + {\widetilde{b^m_n}}^k)\lfloor_{E_i\cap E_j} = (u^m_n, \nabla_{\mu}u^m_n + \widetilde{b^m_n})\lfloor_{E_i \cap E_j}.$
\end{itemize}
All of the above demands are satisfied immediately, because of the smoothness of $u^m_n$ and the choice of the vector field $\widetilde{b^m_n}.$
As the listed conditions are satisfied, by application of the mentioned characterisation, we have $(u^m_n,\widetilde{b^m_n})\in D(A_{\mu}).$

Now we consider the diagonal sequence $$(u^n_n,\widetilde{b^n_n})_{n \in \mathbb{N}}.$$

By estimate \eqref{hestimate}, the construction of the Cosserat field $\widetilde{b^n_n}$ and a completeness of the space $H^1_{\mu}$ it follows that 
$$\lim_{n \to \infty} \left( \nabla_{\mu}u^n_n +\widetilde{b^n_n} \right) \in \left(H^1_{\mu}\right)^3.$$

This means that 
$$A_{\mu}(u^n_n,\widetilde{b^n_n}) \xrightarrow{\left(L^2_{\mu}\right)^{3\times 3}} D \in \left(L^2_{\mu}\right)^{3\times 3}.$$ 
We know that $u^n_n \xrightarrow{L^2_{\mu}} u \in L^2_{\mu}.$ Moreover, by the way it was constructed, $\widetilde{b^n_n}$ converges in $L^2_{\mu}(\R^3;T^{\perp}_{\mu})$ to some element, from now on called as $\widetilde{b},$ 
 that is  $\widetilde{b^n_n} \xrightarrow{L^2_{\mu}(\R^3;T^{\perp}_{\mu})} \widetilde{b}.$ 
A compactness of the operator $A_{\mu}$ (Theorem \ref{main}) implies that
$$A_{\mu}(u^n_n,\widetilde{b^n_n}) \xrightarrow{\left(L^2_{\mu}\right)^{3 \times 3}} A_{\mu}(u,\widetilde{b})$$
and 
$$(u,\widetilde{b}) \in D(A_{\mu}).$$

Finally, under the assumption $L_{\mu}u^n_n \xrightarrow{L^2_{\mu}} B,$ we obtain that 
$$L_{\mu}u^n_n \xrightarrow{L^2_{\mu}} L_{\mu}u,$$ 
thus 
$$L_{\mu}u_n \xrightarrow{L^2_{\mu}} L_{\mu}u$$ 
and $u \in D(L_{\mu}).$

On each component manifold $E_k,\; k=1,2,$ we have shown the convergence ${\norm{u^n_n - u_n}_{H^2(E_k)} \to 0,}$ as $n \to \infty.$ This implies that on each boundary $\partial E_k$ we have a convergence of normal traces in $L^2(\partial E_k).$ In this way we conclude that both
$$(u,\widetilde{b}) \in D(A_{\mu})_N\; \text{ and }\; u \in D(L_{\mu})_N.$$
This proves that the operator $L_{\mu}$ is closed in the case of two components of $\dim E_i = \dim E_j =2.$

$\boldsymbol{\dim E_i=\dim E_j =1}$

Finally, we study the remaining case of $\dim E_i=\dim E_j =1.$

The simplicity of a one-dimensional junction gives chances to describe explicitly the domain $D(A_{\mu})$ (see Example 5.2 in \cite{Bouchitte}). Let $P=E_i\cup E_j,$ $\tau$ be a tangent unit vector field and $\nu$ be a normal outer unit vector field. As it is stated in paper \cite{Bouchitte} in Example 5.2 we have the following characterisation of the space $D(A_{\mu})$ in a case of structures consisting only one-dimensional components:
$$D(A_{\mu}) = \{(u,b)\in H^1_{\mu}\times L^2_{\mu}(\R^3;T^{\perp}_{\mu}): \text{ for } k=\{i,j\}, u\lfloor_{E_k}\in H^2(E_k), b\lfloor_{E_k}\in H^1(E_k), u'\tau + b\nu\in C(P)\}.$$

Without loss of generality, we may assume that $E_i\cap E_j = \{(0,0)\} \in \R^2.$ For $k \in \{i,j\},$ let $\widehat{v}:=\eta_{E_k}v$ denote the parallel transport on the component manifold $E_k$ of a vector $v\in \R^2.$ We put $b_{E_i}(0,0) := (\nabla_{\mu}u)\lfloor_{E_j}$ and we define the Cosserat vector field on $E_i$ as $\widehat{b}_{E_i}:=\eta_{E_i}b_{E_i}(0,0).$ Analogously we proceed on $E_j.$ Then we put $\widehat{b}:= \widehat{b}_{E_i} + \widehat{b}_{E_j}.$ By the above characterisation of $D(A_{\mu})$ it is easy to see that for any element of the sequence $(u_n,b_n) \in D(A_{\mu})$ there exists a mentioned before modification $\widehat{b}_n$ such that $(u_n,\widehat{b}_n)\in D(A_{\mu}).$ By standard properties of Sobolev functions (e.g., see \cite{leo} Theorem 7.13), if $v \in H^1(E_k),$ then $v \in C(E_k),$ and moreover, if $v_n \in H^1(E_k)$ and $v_n \xrightarrow{H^1(E_k)}v,$ then exists a subsequence converging locally uniformly. Each $\widehat{b}_n$ is well-defined and a vector field $b$ being a $L^2_{\mu}$-limit of the sequence of this Cosserat vectors exists. Due to the fact that the $L^2_{\mu}$-convergence of $L_{\mu}u_n$ implies convergence of $u_n$ in the $H^2(E_k)$-sense, and further that $(u,b)\in D(A_{\mu}),$ by closedness of the operator $A_{\mu}$ in $L^2_{\mu}$ we conclude that $u \in D(L_{\mu})$ and $B=L_{\mu}u.$ This proves closedness in the sense of $L^2_{\mu}$-norm convergence of the operator $L_{\mu}.$
The same argument as the one used in the previously considered case gives that 
the result is valid also in the subspace $D(L_{\mu})_N.$

\textbf{``From local to global''}

We have examined closedness in the ``local'' variants. To show a similar result on a whole structure, we need to gather all established results.

Let $u_n \in D(L_{\mu})_N,$ $u_n \xrightarrow{L^2_{\mu}}u,$ $L_{\mu}u_n \xrightarrow{L^2_{\mu}} B.$ Denote $u^p_n:=u_n\lfloor_{\supp \mu_p}.$  We have shown that locally for all $p \in I$ exists $u^p \in D(L_{\mu_p})_N$ which satisfies $u^p_n \xrightarrow{L^2_{\mu_p}}u^p$ with  $L_{\mu_p}u^p_n \xrightarrow{L^2_{\mu_p}} L_{\mu_p}u^p.$ Using the introduced partition of unity, we write $$u_n = \sum_{p \in I}\alpha_p u^p_n \text{ and } u = \sum_{p \in I}\alpha_p u^p.$$ By the aforementioned observations we have for all $p\in I$ 
\begin{itemize}
\item $\alpha_pu^p_n \in D(L_{\mu_p})_N,$
\item $\alpha_p u^p_n \xrightarrow{L^2_{\mu_p}}\alpha_pu^p,$
\item $L_{\mu_p}(\alpha_pu^p_n) \xrightarrow{L^2_{\mu_p}} L_{\mu_p}(\alpha_pu^p).$
\end{itemize}
Applying Proposition \ref{proppp} gives $u \in D(L_{\mu})_N$ satisfying 
$L_{\mu}u_n \xrightarrow{L^2_{\mu}} L_{\mu}u.$ This proves that the operator $L_{\mu}$ is closed.

\qed
\end{dow2}

The next propositions show more properties of the operator $L_{\mu}.$ For the need of our future applications, instead of the operator $L_{\mu}$ we consider the operator $Id - \alpha L_{\mu},$ for some fixed positive constant $\alpha.$ 

\begin{prop}
For any $\alpha >0,$ the operator $Id - \alpha L_{\mu}: D(L_{\mu})_N \to L^2_{\mu}$ is closed with respect to the $L^2_{\mu}$-convergence.
\begin{dow}
A sum of closed operators is a closed operator. The identity operator $Id$ is continuous, thus it is closed, and Theorem \ref{main} provides that the operator $L_{\mu}$ is closed. 
\qed
\end{dow}
\end{prop}

\begin{prop}\label{prop4}
For any $\alpha \in (0,+\infty),$ the operator $Id - \alpha L_{\mu}: D(L_{\mu})_N \to L^2_{\mu}$ is injective. 
\begin{dow}
Proceeding by a contradiction, let us assume that there exists a positive eigenvalue $\lambda = \frac{1}{\alpha}$ of the operator $L_{\mu}$ with the corresponding eigenfunction $w \in D(L_{\mu}).$ Multiplying the eigenvalue equation $\lambda w = L_{\mu} w$ by the eigenfunction $w,$ using the fact, that as $w \in D(L_{\mu})$ the value of $L_{\mu}$ can be computed componentwise, integrating by parts and by the assumption $(B\nabla_{\mu}u\cdot n)\lfloor_{\partial \Omega} = 0$ we obtain 
$$\lambda\norm{w}^2_{L^2_{\mu}} = \int_{\Omega} L_{\mu} w w d\mu = -\int_{\Omega} B\nabla_{\mu}w \nabla_{\mu} w d\mu.$$ 
By ellipticity of the operator $B$ (or equivalently by parabolicity of $L_{\mu}$) we estimate the integral $\int_{\Omega} B\nabla_{\mu}w \nabla_{\mu} w d\mu$ from above and conclude that
$$-\theta \int_{\Omega} |\nabla_{\mu} u|^2 d\mu \geqslant \lambda\norm{w}^2_{L^2_{\mu}}.$$
This implies, that $\lambda < 0,$ what is a contradiction. 
\qed
\end{dow}
\end{prop}

\begin{prop}
$D(L_{\mu})_N$ is dense in $L^2_{\mu}.$
\begin{dow}
We have the inclusion $\{u\in C^{\infty}(\Omega): (B\nabla_{\mu}u\cdot n)\lfloor_{\partial \Omega}=0\} \subset D(L_{\mu})_N.$ Moreover, the space ${\{u\in C^{\infty}(\Omega): (B\nabla_{\mu}u\cdot n)\lfloor_{\partial \Omega}=0\}}$ is a dense subset of $L^2_{\mu}.$ This two facts imply that $D(L_{\mu})_N$ is dense in $L^2_{\mu}.$ 
\qed
\end{dow}
\end{prop}

Now we can proceed to the construction of a semigroup generated by the operator $\Delta_{\mu}.$

\subsection{Generation of semigroup}
Well-known methods of constructing semigroups generated by linear operators are based on the notion of the resolvent operator. In the case considered in this paper, the operator $\lambda Id - \Delta_{\mu}: D(\Delta_{\mu})_N \to L^2_{\mu}$ (or equivalently the operator $Id - \alpha \Delta_{\mu}$) is not surjective, thus its inverse does not exist and it is not possible to introduce the resolvent operator (at least in a standard way). To circumvent this problem we apply a technique proposed in \cite{Mag}, where instead of inverting the given operator, the ``forward'' iterations on a set of ``very regular'' functions are considered.

Throughout this section, we consider the operator $\Delta_{\mu},$ but it seems possible to generalise the presented method to cover the case of a more general form of the operator $L_{\mu}.$  

In further reasonings the space $A_s(\Delta_{\mu})_N$ plays a role of a part of the domain $D(\Delta_{\mu})_N$ in which all needed operations make sense. By this fact the set $A_s(\Delta_{\mu})_N$ will naturally replace the set $D(\Delta_{\mu})_N.$

\begin{defi}\label{veryreg}(Domain of iterations, ``very regular'' functions)\\
Denote $D(\Delta^1_{\mu})_N:= D(\Delta_{\mu})_N.$ For $n \in \{2,3,...\}$ a domain of the $n$-th iteration of the operator $\Delta_{\mu}$ is defined as 
$$D(\Delta^n_{\mu})_N:= \bigl\{u \in L^2_{\mu}: \Delta^{n-1}_{\mu}u \in D(\Delta_{\mu}),\; u \in D(\Delta^{n-1}_{\mu})_N\bigr\}.$$ Put $D_N:= \bigcap_{n=1}^{+\infty} D(\Delta^n_{\mu})_N.$ 

We introduce a space of ``very regular'' functions as 
$$A_s(\Delta_{\mu})_N:=\Bigl\{u \in D_N: \exists M>0\; \forall n \in \mathbb{N}\; \forall t \in (0,s)\; \frac{t^n}{n!}\norm{\Delta_{\mu}^n u}_{L^2_{\mu}}<M\Bigr\}.$$ Moreover, let 
$L^2_{\mu, s}:= \overline{A_s(\Delta_{\mu})_N}^{\norm{\cdot}_{L^2_{\mu}}}.$ 
It will be useful to put $\mathcal{E}(\Delta_{\mu})^T_N:= \bigcap_{s\in [0,T]}A_s(\Delta_{\mu})_N.$
\end{defi}  

We will need to introduce a notion of smoothness with respect to an operator that generates a semigroup. In these terms, the elements of the set $A_s$ defined above can be interpreted as functions analytic with respect to a given differential operator.  

\begin{defi}(Smooth functions with respect to a generator)\\
Let $V$ be a semigroup on the space $L^2_{\mu},$ let $G$ be a generator of $V.$ We define the space $$C^{\infty}(V):=\{u \in \bigcap_{n=1}^{\infty}D(G^n): \norm{G^n u}_{L^2_{\mu}}<+\infty\; \forall n\in \{0,1,2,3,...\}\}.$$ The family of seminorms $\{\norm{G^n \cdot}_{L^2_{\mu}}\}_{n\in \{0,1,2,...\}}$ will be denoted for short as $\norm{\cdot}_{C^{\infty}(V)}.$
In the case of the operator $\Delta_{\mu}$ we introduce the space
$$C^{\infty}(\Delta_{\mu})_N:=\{u \in {D}_N: \forall n \in \{0,1,2,...\},\; \norm{\Delta_{\mu}^n u}_{L^2_{\mu}}<+\infty\}.$$ Similarly, we denote the related family of seminorms by $\norm{\cdot}_{C^{\infty}(\Delta_{\mu})_N}.$
\end{defi}

\begin{prop}\label{ein}
Assume that the operator $G$ generates the semigroup $V,$ $\Delta_{\mu} \subset G$ be a restriction of $A.$ Then for any $u \in A_s(\Delta_{\mu})_N$ we have $\norm{u}_{C^{\infty}(V)}=\norm{u}_{C^{\infty}(\Delta_{\mu})_N}.$
\begin{dow}
An elementary observation. 
\qed
\end{dow} 
\end{prop}

Next proposition shows that spaces $L^2_{\mu,s}$ are independent of the parameter $s>0.$ 

\begin{prop}\label{rownasie}
For any $s>0$ the equality $L^2_{\mu,s} = L^2_{\mu}$ holds.
\begin{dow}
Firstly, let us note that any function $\phi \in C^{\infty}(E_1)$ can be extended to the structure $E=E_1 \cup E_2$ by the formula 
$\widetilde{\phi} =\begin{cases}
\phi, & E_1\\
\phi, & E_2
\end{cases}$ 
is a member of $\bigcap_{n=1}^{\infty}D(\Delta^n_{\mu}).$ This is a result of Proposition 3.11 in \cite{Bouchitte} characterising membership in the set $D(\Delta_{\mu}).$ Choosing a suitable family of smooth functions, we conclude that any function 
$\widetilde{v} =\begin{cases}
v, & E_1\\
v, & E_2
\end{cases},$
where $v \in L^2(E_1),$ satisfies $\widetilde{v} \in \overline{A_s(\Delta_{\mu})_N}^{\norm{\cdot}_{L^2_{\mu}}}.$ To show that the set $\overline{A_s(\Delta_{\mu})_N}^{\norm{\cdot}_{L^2_{\mu}}}$ contains other functions of the class $L^2_{\mu}$ we will follow the following procedure. We indicate a family of smooth functions $\phi_n \in L^2(E_1)$ satisfying Neumann boundary conditions on $E_1.$ We also demand that each $\phi_n$ can be extended by zero to the whole $E,$ that is 
$\widetilde{\phi_n} =\begin{cases}
\phi_n, & E_1\\
0, & E_2
\end{cases}$ and that $\widetilde{\phi_n}\in A_s(\Delta_{\mu})_N.$ Besides that, the family $\{\phi_n\}$ should span a ``large enough'' subspace of $L^2_{\mu},$ for instance, the odd subspace of $L^2_{\mu}.$ Next, we show that characteristic functions
$\widetilde{\chi_A} =\begin{cases}
\chi_A, & E_1\\
0, & E_2
\end{cases}, A \subset E_1$
can be obtained as a $L^2_{\mu}$-limit of elements of $A_s(\Delta_{\mu})_N.$ This is enough to conclude that $L^2_{\mu}=\overline{A_s(\Delta_{\mu})_N}^{\norm{\cdot}_{L^2_{\mu}}}.$ Due to an explicit construction of the sequence $\{\phi_n\}$ we will need to conduct the mentioned procedure separately for the case $\dim E_1=\dim E_2=1$ and the case $\dim E_1=\dim E_2=2.$ Let us start with the one-dimensional situation. We define the family $\{H_{2k-1}, H_{4k-2}, H_{4k}: k=1,2,3,...\},$ where $H_{4k-2}=
\begin{cases}
\sin(\frac{\pi}{2}nx)+ \frac{\pi}{2}nx, & E_1\\
0, & E_2
\end{cases},$
$H_{4k}=
\begin{cases}
\sin(\frac{\pi}{2}nx)- \frac{\pi}{2}nx, & E_1\\
0, & E_2
\end{cases},$
$H_{2k-1}=
\begin{cases}
\sin(\frac{\pi}{2}nx), & E_1\\
0, & E_2
\end{cases}.$\\
Let $\{\alpha_n\}\in l^2$ be a sequence of coefficients of the Fourier series of a fixed odd function $f \in L^2(E_1),$ that is $S_N = \sum_{n=1}^N \alpha_n \sin(\frac{\pi}{2}nx)$ and $\norm{S_N - f}_{L^2} \to 0$ as $N \to \infty.$ Let us assume the decay condition $|\alpha_n|\leqslant \frac{C'}{n^4},$ for some positive constant $C'.$ Denote $h_n := H_n\lfloor_{E_1}$ and we introduce $\gamma_N:= \sum_{n=1}^N \alpha_n h_n.$ There exists a constant $c \in \R$ for which $\norm{\gamma_N-S_n -cx}_{L^2} \to 0$ as $N \to \infty.$ Thus $\gamma_N \xrightarrow[N\to \infty]{L^2} f + cx.$ Our aim is to replace the function $f$ with the function $cx$ and conduct a similar procedure. A series of absolute values of Fourier coefficients of the function $cx$ decays asymptotically to $\frac{1}{n},$ thus the decay condition needed to obtain the above convergence of partial sums is violated, and we cannot proceed directly. As odd periodic smooth functions are dense in $L^2(E_1)$ (assuming that the endpoints $-1$ and $1$ of $E_1$ are glued together, thus we consider smoothness in the sense of the torus $\T^1$) let us take a sequence $\{z_m\}$ of such functions converging in the $L^2$-norm to the function $cx.$ It is clear that the Fourier coefficients $\{\widehat{z_m}(k)\}$ of each $z_m$ decay faster then $\frac{C}{k^p}$ for any exponent $p \in \mathbb{N}.$ Let $\zeta^m_N:=\sum_{k=1}^N \widehat{z_m}(k) h_k.$ We have $\zeta^m_N \xrightarrow[N \to \infty]{L^2} z_m + c_mx.$ By the fact that the $L^2$-convergence of functions implies $l^2$-convergence of corresponding Fourier coefficients we judge that $\zeta_N^N \xrightarrow[N \to \infty]{L^2} cx + \tilde{c}x=(c+\tilde{c})x,$ where $c,\tilde{c}$ are real constants. Now we see that 
$\begin{cases}
(c+\tilde{c})x, & E_1\\
0, & E_2
\end{cases}$ is in $\overline{A_s(\Delta)_N}^{\norm{\cdot}_{L^2_{\mu}}}.$
Without any loss on generality, we might expect that both constants $c$ and $\widetilde{c}$ are non-zero. From the above observations we immediately derive that $\gamma_N - \frac{c}{c+\widetilde{c}}\zeta^N_N \xrightarrow[N \to \infty]{L^2} f,$ so the function 
$\begin{cases}
f, & E_1\\
0, & E_2
\end{cases}$ 
is a member of the space $\overline{A_s(\Delta_{\mu})_N}^{\norm{\cdot}_{L^2_{\mu}}}.$ Repeating this procedure on the component $E_2$ we prove that $L^2_{odd}(E_1) \times L^2_{odd}(E_2) \subset \overline{A_s(\Delta_{\mu})_N}^{\norm{\cdot}_{L^2_{\mu}}},$ where $L^2_{odd}$ stands for the odd subspace of the space $L^2.$ To finish the study of this instance, we notice that both the constant function $1 \in \overline{A_s(\Delta_{\mu})_N}^{\norm{\cdot}_{L^2_{\mu}}}$ and the function 
$\begin{cases}
\begin{cases}
1, & x\in A\\
-1, & x \notin A 
\end{cases}, & E_1\\
0, & E_2
\end{cases} \in \overline{A_s(\Delta_{\mu})_N}^{\norm{\cdot}_{L^2_{\mu}}}$ for any measurable set $A \subset E_1.$ As a consequence, we obtain that for any measurable $A \subset E_1$ we have 
$\begin{cases}
\chi_A, & E_1\\
0, & E_2
\end{cases} \in \overline{A_s(\Delta_{\mu})_N}^{\norm{\cdot}_{L^2_{\mu}}}.$ This gives that $L^2_{\mu} = \overline{A_s(\Delta_{\mu})_N}^{\norm{\cdot}_{L^2_{\mu}}}$ in the one-dimensional case.

In the two-dimensional case, we will follow an essentially similar way. The crucial difficulty is finding a disc analogue of the functions $h_n.$ It turns out that all the necessary requirements are fulfilled by the family of functions $\{u_n, n=1,2,3,...\}$ which in polar coordinates can be written as $u_n(\phi,r) = c^{-1}_n \sin(n\phi)J_n(j_n'r),$ where $c^{-1}_n$ is a $L^2$-normalizing constant, $J_n$ is the n-th order Bessel function of the first kind, and $j_n'$ is a root of $J_n'.$ For a more detailed discussion of the family $u_n$ see \cite{wats}. As $u_n$ is an eigenfunction of the Neumann Laplacian corresponding to the eigenvalue $j_n',$ by the standard facts of the spectral theory we know that $u_n$ are smooth, satisfy Neumann boundary conditions and the functions $u_n$ span the odd subspace of $L^2(E_1).$ Moreover, the eigenvalue equality implies that all derivatives of each $u_n$ have proper growth. The polar representation of $u_n$ shows that each $u_n$ and each $\Delta^m u_n$ vanishes on the intersection set $\Sigma = E_1 \cap E_2$ (or this happens up to rotation), thus $u_n$ can be extended by zero to the whole structure $E$ and 
$\widetilde{u}_n=\begin{cases}
u_n, & E_1\\
0, & E_2
\end{cases}$
satisfies $\widetilde{u}_n \in A_s(\Delta_{\mu})_N.$ Now, we are able to proceed analogously to the earlier discussed one-dimensional case. In this way we conclude $L^2_{\mu} = \overline{A_s(\Delta_{\mu})_N}^{\norm{\cdot}_{L^2_{\mu}}}$ what finishes the proof. 
  
\qed
\end{dow}
\end{prop}	
	
The fact that convergence in the sense of the graph of $\Delta_{\mu}$ implies convergence in the Sobolev norm $H^2$ on each component manifold will be needed in further considerations.

\begin{prop}\label{zwei}
Assume that $u_n \in A_s(\Delta_{\mu})_N,$ both the sequence $u_n$ and the sequence $\Delta_{\mu} u_n$ are convergent in the $L^2_{\mu}$-norm, then the sequence $u_n$ is convergent in $H^2(E_i),\; i=1,2.$
\begin{dow}
The evoked earlier characterisation (see Proposition 3.11 in \cite{Bouchitte}) of the space $D(\Delta_{\mu})$ shows that membership $u_n \in A_s(\Delta_{\mu})_N$ implies $u_n \in H^2(E_i),$ for $i=1,2.$ By considering separately each component $E_i,\; i=1,2$ and using the local regularity estimates (see for instance \cite{Evans}, p.306) we obtain that convergence in the sense of graph of $\Delta$ implies convergence in the seminorm ${\norm{\nabla^2 \cdot}_{L^2_{\mu}}.}$By the interpolation theorem (see thm 5.2, p.133 in \cite{adams}) we obtain that a joint convergence in seminorms $\norm{\cdot}_{L^2(E_i)}$ and $\norm{\nabla^2\cdot}_{L^2(E_i)}$ is equivalent to a convergence in the standard $H^2(E_i)$-norm. 
\qed
\end{dow}
\end{prop}	
	
To provide that a solution constructed by the action of a semigroup is fully valuable, we need to ensure that for each time $t$ it is a member of the domain of the generator. 

\begin{lem}
Let $u \in A_s(\Delta_{\mu})_N,$ let $V_s$ be a family of operators constructed in the proof of Theorem 2 in \cite{Mag}. Then for each $t \in (0,s)$ we have $V_s(t)u \in D(\Delta_{\mu})_N.$
\begin{dow}
From the proof of Theorem 2 in \cite{Mag} (p. 99) it follows that for an arbitrary $u \in A_s(\Delta_{\mu})_N$ and $t \in (0,s)$ we have $V_s(t)u \in \overline{A_s(\Delta_{\mu})_N}^{\norm{\cdot}_{C^{\infty}(V_s)}}.$ By the equality of Proposition \ref{ein} it follows $V_s(t)u \in \overline{A_s(\Delta_{\mu})_N}^{\norm{\cdot}_{C^{\infty}(\Delta_{\mu})_N}}.$ Let $w_n \in A_s(\Delta_{\mu})_N$ be a Cauchy sequence in $\norm{\cdot}_{C^{\infty}(\Delta_{\mu})_N}.$ As Proposition \ref{zwei} implies convergence of $w_n$ in the sense of $H^2(E_i)$ on each component manifold $E_i,\; i=1,2,$ and as the normal trace $\frac{\partial w_n}{\partial n}\lfloor_{\partial E_i}:H^1(E_i) \to L^2(E_i), i=1,2$ is continuous we conclude that the limit $w:=\lim_{n \to \infty}w_n$ in $\norm{\cdot}_{C^{\infty}(\Delta_{\mu})_N}$ satisfies 
$\frac{\partial w}{\partial n}\lfloor_{\partial E_i}=0,\; i=1,2.$ It is easy to check that $w \in H^2(E_i)$ and $w \in C(E).$ Now using the reasoning presented in the proof of closedness of the operator $L_{\mu}$ we deduce that $w \in D(\Delta_{\mu})_N.$ In this way we obtained $V_s(t)u \in D(\Delta_{\mu})_N,$ for all $t \in (0,s).$
\qed
\end{dow}
\end{lem} 	

We will also need the following simple observation.

\begin{prop}
$\mathcal{E}(\Delta_{\mu})^T_N$ is dense in $L^2_{\mu}.$
\begin{dow}
It follows directly from the proof of Proposition \ref{rownasie}.
\qed
\end{dow}
\end{prop}	

The following proposition is equivalent to the continuity of the resolvent operator if such an operator is well-defined.
	
\begin{prop}\label{drei}
For any neighbourhood of zero $W \subset L^2_{\mu}$ exists a neighbourhood of zero $U \subset L^2_{\mu}$ such that if for all $k\in \{1,2,...\},$ for all $u \in D(\Delta_{\mu}^k)_N$ and some constants $\alpha,C>0$ we have $(1-\alpha C)^{-k}(Id-\alpha \Delta_{\mu})^k u\in U,$ then $u \in W.$
\begin{dow}
Let us recall that the operator $Id - \alpha \Delta_{\mu}:D(\Delta_{\mu})_N \subset L^2_{\mu} \to L^2_{\mu}$ is closed due to closedness of $\Delta_{\mu}.$ Assume that the statement postulated in the thesis is not true. This means that the following sentence is valid:
$\exists k \in \{1,2,...\}\; \exists u \in D(\Delta_{\mu}^k)_N,$ exists a neighbourhood of zero $W \subset L^2_{\mu}$ such that for any neighbourhood of zero $U \subset L^2_{\mu}$ we have $(1-\alpha C)^{-k}(Id-\alpha \Delta_{\mu})^k u \in U \wedge u \notin W.$ This implies existence of a sequence $u_n \in D(\Delta_{\mu})_N$ such that $\norm{(Id-\alpha \Delta_{\mu})u_n}_{L^2_{\mu}} \to 0$ and $\norm{u_n}_{L^2_{\mu}}>c'>0.$ By continuity of the classical resolvent operator defined for the Laplace operator $\Delta$ on each component $E_i,$ we see that the convergence $\norm{(Id-\alpha \Delta_{\mu})u_n}_{L^2_{\mu}} \to 0$ implies $\norm{u_n}_{L^2_{\mu}} \to 0$ as $n \to \infty.$
By induction over $k \in \{1,2,...\},$ we stipulate that the condition expressed in the thesis is valid. 
\qed
\end{dow}
\end{prop}	
	
Finally, we can summarise our study over generating semigroups in the following theorem.

\begin{tw}\label{eqisemi}
The operator $\Delta_{\mu}$ generates equicontinuous semigroup $S:[0,T] \times L^2_{\mu} \to L^2_{\mu}.$
\begin{dow}
We verified, that on the interval $[0,T]$ all conditions of Theorem 2 of \cite{Mag} (see p.95) are satisfied. This implies, that $\Delta_{\mu}$ generates an equicontinuous semigroup $S:[0,T] \times L^2_{\mu} \to L^2_{\mu}.$
\qed
\end{dow}
\end{tw}

We may conclude the results established in this section with the proof of Theorem \ref{existence}.

\begin{dow1}
Let $S$ be the semigroup as in Theorem \ref{eqisemi}. As $g \in A_T(\Delta_{\mu})_N$ and the semigroup $S$ is generated by the operator $\Delta_{\mu}$ we immediately notice that the function $u(t,x):=(S(t)g)(x)$ is a solution in a sense of Definition \ref{parabolic2}.
\qed
\end{dow1}

\begin{przyk}
We present an example of an application of Theorem \ref{existence} to the existence of solutions in a case of two orthogonal discs intersecting each other.

Let $\mu \in \mathcal{S}$ be a low-dimensional structure of the form $$\mu:=\mathcal{H}^2|_{D_1}+\mathcal{H}^2|_{D_2},$$ where $D_1:=\{(x,y,0)\in \R^3: x^2+y^2\leqslant 1\}$ and $D_2:=\{(x,0,z)\in \R^3: x^2+y^2 \leqslant 1\}.$ We denote $D:= D_1 \cup D_2$ and $\partial D := \partial D_1 \cup \partial D_2.$

Assume that the matrix $B=(b_{ij}),\; i,j\in\{1,2,3\}$ from Definition \ref{secondop} consists of constant entries $b_{ii}\equiv 1,\; i \in \{1,2,3\}$ and $b_{ij}\equiv 0,\; i\neq j,\; i,j\in\{1,2,3\}.$ The operator $L_{\mu}:D(L_{\mu})_N \to L^2_{\mu}$ related with this matrix is of the simple form $$L_{\mu}u=\sum_{i=1}^3(\nabla^2_{\mu}u)_{ii} = \Delta_{\mu}u.$$

We consider the evolutionary heat equation
\begin{equation}\label{parabolic3}
\begin{aligned}
u_t - \Delta_{\mu}u = 0 \;\;\;(\mathcal{H}^2|_{D_1}+\mathcal{H}^2|_{D_2})\times l^1([0,T])-\text{a.e. in } (D_1 \cup D_2) \times [0,T] \\
(\nabla_{\mu} u,\eta) = 0 \;\text{ on } (\partial D_1 \cup \partial D_2) \times [0,T] \\
u = g \;\text{ on } (D_1 \cup D_2) \times \{0\}.
\end{aligned}
\end{equation}
Here $g \in L^2_{\mu}$ and $\eta$ is a vector field such that $\eta\lfloor_{D_1}$ is the outer normal unit vector field to $D_1$ and $\eta\lfloor_{D_2}$ is the outer normal unit vector field to $D_2.$

By the result of Theorem \ref{existence}, there exists a unique function $u:[0,T]\to L^2_{\mu}$ which is a solution of the heat transfer issue \eqref{parabolic3} in the sense of Definition \ref{parabolic2}.  
\end{przyk}

\section{Non-semigroup approach}	
	A fundamental drawback of the presented semigroup method is that it provides the existence of solutions in a narrow space of functions and under the assumption of a very regular initial input (see Theorem \ref{eqisemi}). This is closely related to two aspects: a form of the considered second-order operator and the non-existence of the corresponding resolvent operator. The second-order operator $L_{\mu}$ was constructed to be consistent with the variational theory introduced in \cite{Bouchitte} and the operators used there. The operator $L_{\mu}$ needs restrictive conditions to be posed on function spaces to provide well-definiteness of it. Further restrictions are caused by the fact that the considered operator is non-invertible, implying that the resolvent operator does not exist (at least in a classical sense, see Subsection 3.1). Due to this, our method of constructing solutions bases on ``forward'' iterations of the considered operator. Such a method demands further significant restrictions on the considered space of functions. 
	Now our goal is to touch the low-dimensional parabolic issues differently. Firstly ensuring that the class of solutions is wide enough, and later on, we try to deduce additional regularity of obtained solutions.
	In this section, we focus on constructing solutions to weak versions of parabolic problems with initial data of low regularity, and in further considerations, we show in what sense the regularity of it can be upgraded.
	We begin by introducing a notion of weak parabolic problems and by examining the existence of solutions in the first-order framework of \cite{Rybka}. \\

 Throughout this section, we can loosen our restrictions on the class of considered low-dimensional structures. We need only to assume that $\mu \in \widetilde{\mathcal{S}}.$ This means that components of low-dimensional structures may have non-fixed dimensions.\\
 
We adapt the framework exposed, for instance, in book \cite{None} to provide the existence of solutions to the weak parabolic problems on the low-dimensional structures. \\

We will need the next proposition to prove the existence of solutions to \eqref{parabolicweak}.

\begin{prop}\label{conde}
For any $\phi \in \mathring{\mathcal{T}}$ we have
\begin{itemize}
\item[a)] $\norm{\phi}_{L^2H^1_{\mu}} \leqslant \norm{\phi}_{{\mathcal{T}}},$
\item[b)] $E(\phi,\phi)\geqslant C\norm{\phi}^2_{{\mathcal{T}}},$ where $C>0$ is a constant independent of a choice of $\phi.$ 
\end{itemize}
\begin{dow}
Point a) follows directly from the definitions of the related norms. This means that the space $\mathring{\mathcal{T}}$ embeds continuously in the space $\mathring{\mathcal{H}}.$ We move to the proof of point b). Using integration by parts with respect to the time variable and ellipticity of the matrix operator $B,$ we see that 
$$E(\phi,\phi) = \int_0^T\int_{\Omega} (B\nabla_{\mu} \phi, \nabla_{\mu} \phi) d\mu dt - \int_0^T \int_{\Omega} \phi \phi' d\mu dt \geqslant C'\left(\norm{\nabla_{\mu} \phi}^2_{L^2L^2_{\mu}} + \norm{\phi(0)}^2_{L^2_{\mu}}\right)$$ with a constant $C'>0$ that does not depend on $\phi.$ Rewriting the $\mathcal{T}$-norm, we estimate
$$\norm{\phi}^2_{\mathcal{T}} = \norm{\phi}^2_{L^2H^1_{\mu}}+\norm{\phi(0)}^2_{L^2_{\mu}}=\norm{\nabla_{\mu}\phi}^2_{L^2L^2_{\mu}}+\norm{\phi}^2_{L^2L^2_{\mu}}+\norm{\phi(0)}^2_{L^2_{\mu}}\leqslant (1+C_p)\norm{\nabla_{\mu}\phi}^2_{L^2L^2_{\mu}}+\norm{\phi(0)}^2_{L^2_{\mu}},$$ where the constant $C_p$ comes from the generalized Poincar{\'e} inequality (formula \eqref{weakpoincare} in Section 2.2). As a conclusion we derive that $C\norm{\phi}^2_{\mathcal{T}}\leqslant E(\phi,\phi).$
\qed
\end{dow}
\end{prop}	

Let us recall the Lions variant of the Lax-Milgram Lemma (see, for instance, \cite{None}).

\begin{tw}\label{lions}
Let $\mathcal{M}$ be a Hilbert space equipped with the norm $\norm{\cdot}_{\mathcal{M}}$ and $\mathcal{N}$ with the norm $\norm{\cdot}_{\mathcal{N}}$ be a normed space. Let $H:\mathcal{M} \times \mathcal{N} \to \R$ be a bilinear form and assume that for any $\phi \in \mathcal{N}$ we have $H(\cdot,\phi)\in \mathcal{M}^*.$ Then the condition $$\inf_{\norm{\phi}_{\mathcal{N}}=1}\sup_{\norm{u}_{\mathcal{M}}\leqslant 1}|H(u,\phi)|\geqslant c >0$$ is equivalent to the fact that for any $F \in \mathcal{N}^*$ there exists $u \in \mathcal{M}$ such that for all $\phi \in \mathcal{N}$ we have $H(u,\phi)=F(\phi).$
\end{tw}

The following proposition (see \cite{None}) serves as a useful criterion for verifying that one of the implications of the Lions Theorem is valid.

\begin{prop}\label{lionscol}
Let there exists a continuous embedding of $\mathcal{N}$ in $\mathcal{M}.$ If there is some positive constant $A$ such that $H(\phi,\phi)\geqslant A\norm{\phi}^2_{\mathcal{N}}$ for all $\phi \in \mathcal{N},$ then for any $F \in \mathcal{N}^*$ exists $u \in \mathcal{M}$ satisfying $H(u,\phi)=F(\phi)$ for all $\phi \in \mathcal{N}.$
\end{prop} 

Now we are prepared to deal with the well-posedness of problem \eqref{parabolicweak}.

\begin{tw}\label{weakpara}
Let $f \in L^2(0,T;{L^2_{\mu}})$ and $u_0 \in \mathring{L^2_{\mu}}.$ The problem \eqref{parabolicweak} has a unique solution $u \in \mathring{\mathcal{H}}.$ 
\begin{dow}
Proposition \ref{conde} implies that the assumption of Proposition \ref{lionscol} is satisfied. Thus by applying the Lions Theorem (Theorem \ref{lions}) the existential part is done. To derive the uniqueness of solutions notice that point b) of Proposition \ref{conde} is stricter than the monotonicity of the operator. Thus the uniqueness is provided by applying Proposition 2.3 of \cite{None} (p. 112).
\qed
\end{dow}
\end{tw}	

As the existence of weak solutions is already discussed, our next goal is to analyse their regularity. 

\subsection{More regular solutions}	

Let us begin with the following, at this moment only formal, computations:
$\int_{\Omega} \nabla_{\mu} u \cdot \nabla_{\mu}v d\mu = \int_{E_1}\nabla_{\mu} u \cdot \nabla_{\mu}v d\bar{x} + \int_{E_2}\nabla_{\mu} u \cdot \nabla_{\mu}v d\bar{x} = 
-\int_{E_1} \Delta u v d\bar{x} + \int_{\partial E_1} \frac{\partial u}{\partial n} v d\sigma - \int_{E_2} \Delta u v d\bar{x} + \int_{\partial E_2} \frac{\partial u}{\partial n} v d\sigma = - \int_{\Omega} Au v d\mu + \int_{\partial E} \frac{\partial u}{\partial n} v d\sigma,$ here we denote $A:= \Delta_{E_1}+ \Delta_{E_2}.$
It should be clear that the above computations make sense only when we additionally assume ${u\in H^2(E_i),\; i=1,2.}$ Moreover, introducing the zero Neumann boundary condition on each component $E_i,\; i=1,2$ it follows by continuity of the operator $A:H^2(E_1)\times H^2(E_2) \to L^2_{\mu}$ that $\int_{\Omega} \nabla_{\mu} u \nabla_{\mu}v d\mu \leqslant C \norm{v}_{L^2_{\mu}},$ for all $v \in H^1_{\mu}.$  

\begin{defi}\label{strop}
Let $D(\mathcal{A}) \subset \mathring{H^1_{\mu}}$ such that for $u \in D(\mathcal{A})$ we have ${\int_{\Omega} (B\nabla_{\mu} u, \nabla_{\mu} v) d\mu \leqslant C \norm{v}_{L^2_{\mu}}}$ for all $v \in \mathring{H^1_{\mu}}$ and some positive constant $C.$
We define the operator ${\mathcal{A}:D(\mathcal{A})\to L^2_{\mu}}$ by the equality $\int_{\Omega} (B\nabla_{\mu} u, \nabla_{\mu} v) d\mu = \int_{\Omega} (\mathcal{A}u) v d\mu$ satisfied for all $v \in \mathring{H^1_{\mu}}.$ Existence of the element $\mathcal{A}u$ is guaranteed by the Riesz Representation Theorem for functionals on a Hilbert space.
\end{defi}	

\begin{kom}
Let us point out that the operator $\mathcal{A}$ is consistent with the integration by parts formula. We do not assume enough regularity on $u$ to ensure the existence of the normal trace in the classical sense. We must follow the generalized approach by applying the results established in \cite{Chen}. This theory provides the existence of the normal trace and the operator $\mathcal{A}$ in the sense of distributions or measures. Further, by a condition from the definition of the domain $D(\mathcal{A})$ (Definition \ref{strop}), we represent the normal trace and other terms of the integration by parts formula (we refer to Theorem 2.2 in \cite{Chen}) as functionals on adequate variants of $L^2$ spaces. In this way we obtain $\int_{\Omega} \mathcal{A}u v d\mu + \int_{\partial \Omega}[B\nabla_{\mu}u,\nu]\lfloor_{\partial \Omega} v d\sigma \leqslant C \norm{v}_{L^2_{\mu}}.$ In a case if $[B\nabla_{\mu}u,\nu]\lfloor_{\partial \Omega}$ is non-zero somewhere on $\partial \Omega,$ then we can find a sequence of continuous functions $v_n \in \mathring{H^1_{\mu}}$ such that the second integral on the left side of the inequality diverges and the other terms on both sides are bounded. This indicates that in fact the operator $\mathcal{A}$ satisfies $\int_{\Omega} B\nabla_{\mu} u \nabla_{\mu} v d\mu = \int_{\Omega} \mathcal{A}u v d\mu.$

Finally, let us observe that ${\{u \in \mathring{H^1_{\mu}}: u|_{E_i}\in H^2(E_i),\; [B\nabla_{\mu}u,\nu]\lfloor_{\partial \Omega}=0,\; i=1,2\} \subset D(\mathcal{A}).}$ 
\end{kom}

\begin{kom}\label{self}
The operator $\mathcal{A}:D(\mathcal{A}) \to L^2_{\mu}$ is selfadjoint. This follows directly from the symmetry of the matrix operator $A_{\mu}(\bar{x})$ for $\mu$-a.e. $\bar{x}\in E$ and the symmetry of the bilinear form $(\cdot,\cdot)_{L^2_{\mu}}.$
\end{kom}

Using the theorem of Lions - Theorem \ref{lions}, we prove the existence of strong solutions in the sense of the operator $\mathcal{A}.$

\begin{prop}\label{exiw}
Let $f \in L^2(0,T;{L^2_{\mu}}), u_0 \in \mathring{H^1_{\mu}}.$ Then there exists a unique $u \in W^{1,2}(0,T;{L^2_{\mu}})$ satisfying $u' + \mathcal{A}u = f$ in $L^2(0,T;L^2_{\mu})$ with $u(0)=u_0.$ Moreover, for almost every $t \in (0,T)$ it is $u(t) \in D(\mathcal{A}).$
\begin{dow}
By a fact that the operator $\mathcal{A}$ is obtained from the non-negative, symmetric bilinear form, see Definition \ref{strop}, and due to its self-adjointness (see Comment \ref{self} above) we are in a case covered by Proposition 2.5 of \cite{None}. Proposition 2.5 of \cite{None} allows adapting the Lions Theorem to the equation with the operator $\mathcal{A}.$ In this way we obtain existence of unique function $u$ such that equality $u' + \mathcal{A}u = f$ is satisfied in the norm of the space $L^2(0,T;L^2_{\mu}).$
\qed
\end{dow} 
\end{prop}

Now we focus on extending a class of accessible initial data $u_0$ and discuss remarks dealing with the regularity of the obtained solution.

\begin{prop}
A solution determined in Proposition \ref{exiw} is a member of the space $C([0,T];\mathring{H^1_{\mu}}).$
\begin{dow}
This is a result of an equivalence of $\int_{\Omega}(B\nabla_{\mu}u,\nabla_{\mu}v)d\mu$ and the scalar product of the Hilbert space $\mathring{H^1_{\mu}}.$ Such equivalence is valid because we assume that on the low-dimensional structure $\mu$ the generalized Poincar{\'e} inequality is true. Moreover, the mentioned equivalence implies that the assumption $u_0 \in \mathring{H^1_{\mu}}$ cannot be relaxed in the given class of solutions to the considered problem.
\qed
\end{dow}
\end{prop}

Finally, we examine the regularity of solutions in the setting where as an initial data we can take the function of the class $\mathring{L^2_{\mu}}.$ 

\begin{prop}
Let the operator $\mathcal{A}$ and the function $f$ be as in Proposition \ref{exiw}. Assume that ${u_0 \in \mathring{L^2_{\mu}}.}$ Then the unique solution $u \in W^{1,2}(0,T;{L^2_{\mu}})$ of the problem
$u'+\mathcal{A}u=f$ in $L^2(0,T;H^{1}_{\mu})$ with $u(0)=u_0 \in \mathring{L^2_{\mu}}$ satisfies $t^{1/2}u' \in L^2(0,T;L^2_{\mu})$ and for an arbitrary $d\in(0,T)$ satisfies $u \in W^{1,2}(d,T;L^2_{\mu}) \cap C([d,T];\mathring{H^1_{\mu}}).$
\begin{dow}
Under the given assumptions, we are in the regime of Corollary 2.4 from \cite{None}, which gives the proposed statements.
\qed
\end{dow}
\end{prop}

A relation between solutions of the form evoked in Proposition \ref{exiw} and weaker solutions recalled in Theorem \ref{weakpara} is investigated in the below lemma.

\begin{lem}
Let $u\in W^{1,2}(0,T;{L^2_{\mu}}), u_0\in \mathring{{H}^1_{\mu}}$ satisfies the equation $u'+\mathcal{A}u=f$ in $L^2(0,T;L^2_{\mu})$ with $u(0)=u_0.$ Then $u$ is a weak solution in the sense of Definition \ref{weak3}. 
\begin{dow}
After multiplying the equation $u'+\mathcal{A}u=f$ by a test function $\phi \in \mathring{\mathcal{C}^{\infty}_0}$ and integrating over $\Omega$ and time, we obtain
$\int_0^T \int_{\Omega}u' \phi + \mathcal{A}u \phi d\mu \dt = \int_0^T \int_{\Omega}u' \phi + (B\nabla_{\mu}u, \nabla_{\mu}\phi) d\mu \dt.$ Continuing by integrating by parts the term with the time derivative it follows that\\
$\int_0^T \int_{\Omega}u' \phi + (B\nabla_{\mu}u, \nabla_{\mu}\phi) d\mu \dt=\int_0^T \int_{\Omega} (B\nabla_{\mu}u, \nabla_{\mu}\phi) - u\phi' d\mu \dt+ \int_{\Omega} u_0 \phi(0)d\mu.$ From this computations we conclude that $E(u,\phi) = F(\phi),$ thus the function $u$ is a weak solution in a sense of Definition \ref{weak3}.
\qed
\end{dow}
\end{lem}

\subsection{Examples}	
The following examples are evoked to present that in the low-dimensional setting, even the simplest stationary elliptic equation might possess weak solutions which give a new quality. We will show that there are solutions that are not simple gluings of classical solutions of projected problems on component manifolds. 

\begin{przyk}
	Assume that $\Omega \subset \mathbb{R}^2$ is a 2-dimensional unit ball $B(0,1),$ $E_1:=\{(y,0): y\in [-1,1]\},$ $E_2:=\{(0,z):z \in [-1,1]\}$ and $\mu := \mathcal{H}^1|_{E_1}+\mathcal{H}^1|_{E_2}.$ Let 
	$f:=\begin{cases}
	y \text{ on } E_1,\\
	0 \text{ on } E_2
	\end{cases}.$  
	Clearly we have $f \in \mathring{L^2_{\mu}}.$	Let us consider the stationary heat problem (see \cite{Rybka} for the existence and uniqueness result)
	\begin{equation}\label{slabe}
	\int_{\Omega} \nabla_{\mu}u\cdot \nabla_{\mu} \phi d\mu = \int_{\Omega} f \phi d\mu
	\end{equation}
	for $\phi \in C^{\infty}_c(\R^2).$
	In this case the projected gradient $\nabla_{\mu}$ has a form of the classical 1-dimensional derivative $\partial$ in variables corresponding to each component $E_i, i=1,2.$ We can write down equation \eqref{slabe} in a form
	\begin{equation*}
	\int_{E_1}\partial_y u \partial_y \phi dy + \int_{E_2}\partial_z u \partial_z \phi dz = \int_{E_1} y \phi dy + \int_{E_2} 0 \phi dz = \int_{E_1} y \phi dy.
	\end{equation*}
	On the component $E_1$ let us take $u_1(y):= -\frac{21}{1080}-\frac{y^4}{12}+\frac{y^3}{6}+\frac{y^2}{6}-\frac{y}{2}.$ Computations show that $u_1'(-1)=u_1'(1)=0,$ $u_1(0)=\frac{-21}{1080},$ and $\int_{E_1}u(y)dy = \frac{7}{180}.$ Let us put on the component $E_2$ the constant function $u_2(z):=-\frac{21}{1080}.$ The function $u_2$ has a mean over $E_2$ equal to $-\frac{7}{180}.$ It is easy to observe that neither $u_1$ nor $u_2$ satisfies the weak equation separately on components, that is $u_i, i=1,2$ is not the solution of $\int_{E_i} \partial_x u_i \partial_x \phi dx = \int_{E_i} f \phi dx,$ for $x\in \{y,z\},$ because $\int_{E_i} u_i dx \neq 0.$ On the other hand, the function 
	$u := \begin{cases}
	u_1 \text{ on } E_1,\\
	u_2 \text{ on } E_2
	\end{cases}$
	belongs to $\mathring{H^1_{\mu}}$ and is a solution to weak problem \eqref{slabe}. Moreover, in the considered case, the uniqueness of solutions to \eqref{slabe} is guaranteed by Theorem 1.1 of \cite{Rybka}. This is equivalent to satisfying that a solution is globally continuous and has a mean zero. 
\end{przyk}

The second example shows that a solution to the weak low-dimensional Poisson problem \eqref{slabe} may differ from a sum of solutions to classical problems considered separately on each component manifold.

\begin{przyk}
	We are considering Poisson equation \eqref{slabe} taking for $\Omega \subset \R^3$ a 3-dimensional unit ball $B(0,1),$ the component manifolds $E_1:=\{(x,y,0):x^2+y^2\leqslant 1\},$ $E_2:=\{(x,0,z): x^2+z^2\leqslant 1\}$ and the measure $\mu := \mathcal{H}^2|_{E_1}+\mathcal{H}^2|_{E_2}.$ Let us define a function $w:E_1 \to \R,$ $w(x,y):=\cos(\pi(x^2+y^2))$ and for the force term let us put 
	$f:= \begin{cases}
	\Delta w \text{ on } E_1,\\
	0 \text{ on } E_2
	\end{cases}.$  We have $\nabla w = -2\pi\sin(\pi(x^2+y^2))(x,y),$ $\frac{\partial w}{\partial n}\lfloor_{\partial E_1} = 0,$ $\Delta w = -4\pi \sin(\pi (x^2+y^2)) - 4 \pi^2 (x^2+y^2) \cos(\pi (x^2+y^2)).$ By a change of variables for the polar coordinates, we can easily verify that $\int_{E_1} \Delta w d\bar{x} =0.$
	The previous observations show that $f \in \mathring{L^2_{\mu}}.$ Both $w$ and the zero function are the solutions of the classical weak Poisson problem posed on the corresponding component, but the function 
	$u:=\begin{cases}
	w \text{ on } E_1\\
	0 \text{ on } E_2
	\end{cases}$
	does not belong to $\mathring{H^1_{\mu}}$ or even to $H^1_{\mu}.$ It should be noted that even the modified function 
	$\widetilde{u}:=\begin{cases}
	w +c_1 \text{ on } E_1\\
	0 +c_2 \text{ on } E_2
	\end{cases}$
	for any $c_1, c_2 \in \R$ is not a member of $H^1_{\mu}.$ The general existence theorem established in \cite{Rybka} can be applied to this problem and provides the existence of the unique solution $\widehat{u}\in \mathring{H^1_{\mu}}.$ The conducted reasoning shows that the solution $\widehat{u}$ will be different from any function $\widetilde{u}$ and thus will be different from the function $u$ in an essential way, meaning that neither only by some constant nor by two different constants added independently on components. 
\end{przyk}

\subsection{Asymptotic convergence}
We focus on studying the asymptotic behaviour of weak solutions to parabolic problems. Using the existential theorems established in Section 4.1, the variational methods of \cite{Rybka} and transferring the results of \cite{goldie}, we prove that solutions to parabolic problem \eqref{parabolicweak} converge as $t \to \infty$ to a solution of the stationary heat equation. 

We will work with a slightly different formulation of the weak equation. Thus we will need an equivalence (see Proposition 2.1 in \cite{None}) phrased by the below proposition. 

\begin{prop}\label{edef}
Weak formulation of the parabolic initial problem \eqref{parabolicweak} is equivalent to a problem of finding $u \in \mathring{\mathcal{H}}$ with $u(0)=u_0 \in L^2_{\mu}$ such that the equation $u'(t)+\mathcal{A}u(t)=f(t)$ is satisfied in the sense of $\mathring{\mathcal{H}}^*$ for a.e. $t\in (0,+\infty).$
\begin{dow}
This can be proved by applying Proposition 2.1 of \cite{None} to equation \eqref{parabolicweak} in the corresponding framework.
\qed
\end{dow}
\end{prop} 

To shorten the notation, we introduce a notion of the energy functional related to the considered problem.

\begin{defi}
Let a functional $E_{\mu}:\mathring{H^1_{\mu}} \to \R$ be defined as $E_{\mu}(u):=\int_{\Omega} (B\nabla_{\mu}u,\nabla_{\mu}u)d\mu.$ The functional $E_{\mu}$ is called the energy functional.
\end{defi} 

To avoid confusion, we formulate the notion of a weak elliptic problem, which we consider further.

\begin{defi}\label{stat}
Let $f \in \mathring{H^1_{\mu}}.$ We say that a function $u^* \in \mathring{H^1_{\mu}}$ is a solution of the stationary heat equation if 
$$\int_{\Omega} (B\nabla_{\mu}u^*,\nabla_{\mu}v)d\mu=\int_{\Omega} fv d\mu$$
is satisfied for all functions $v \in H^1_{\mu}.$
\end{defi}

Adapting the proof of Theorem 1 of \cite{goldie}, (p.269) to the considered case, we obtain the following theorem dealing with the asymptotic behaviour of parabolic solutions. 

\begin{tw}
Let $u^* \in \mathring{H^1_{\mu}}$ be a solution to stationary heat equation (Definition \ref{stat}). A solution $u$ to parabolic issue \eqref{parabolicweak} converges as $t \to \infty$ to $u^*$ in the sense of the $H^1_{\mu}$-norm.
\begin{dow}
A main point of the proof is in verifying that $E_{\mu}(u(t)-u^*) \to 0$ as $t \to \infty.$ Let us compute: 
\begin{equation}\label{differ}
\begin{aligned}
\frac{d}{dt}\norm{u-u^*}^2_{L^2_{\mu}}=\int_{\Omega}u_t(u-u^*)d\mu=\int_{\Omega}f(u-u^*)d\mu-\int_{\Omega}(B\nabla_{\mu}u,\nabla_{\mu}(u-u^*))d\mu=\\ \int_{\Omega}f(u-u^*)d\mu-\int_{\Omega}(B\nabla_{\mu}u,\nabla_{\mu}(u-u^*))d\mu+\int_{\Omega}(B\nabla_{\mu}u^*,\nabla_{\mu}(u-u^*))d\mu-\int_{\Omega}f(u-u^*)d\mu=\\-\int_{\Omega}(B\nabla_{\mu}u,\nabla_{\mu}(u-u^*))d\mu+\int_{\Omega}(B\nabla_{\mu}u^*,\nabla_{\mu}(u-u^*))d\mu =\\ -\int_{\Omega}(B\nabla_{\mu}(u-u^*),\nabla_{\mu}(u-u^*))d\mu=-E_{\mu}(u-u^*)\leqslant 0.
\end{aligned}
\end{equation}
This shows that the term $\norm{u-u^*}^2_{L^2_{\mu}}$ is non-increasing. Please note that if $E_{\mu}(w)=0$ for some $w\in \mathring{H^1_{\mu}},$ then $w \equiv 0.$ This means if $u(t) = u^*$ for some $t \in (0,\infty),$ then this implies we have $u(t+s)=u_0$ for all $s>0.$ Now we see that as $\norm{u-u^*}^2_{L^2_{\mu}}\geqslant 0,$ there exists a sequence of times $t_i \in (0,\infty)$ such that $\frac{d}{dt}\norm{u_i-u^*}^2_{L^2_{\mu}} \to 0$ as $t_i \to \infty.$ Here we denote $u_i:=u(t_i).$
In this way we get $E_{\mu}(u_i-u^*) \to 0$ as $t_i\to \infty.$ As \eqref{differ} shows, the term $\norm{u-u^*}^2_{L^2_{\mu}}$ is non-increasing with respect to $t\in (0,\infty)$ thus we can deduce that the convergence occurs for any sequence of times, that is $E_{\mu}(u-u^*) \to 0$ as $t\to \infty.$ The ellipticity of the operator $B$ implies that $\nabla_{\mu}u \xrightarrow{t \to \infty} \nabla_{\mu}u^*$ in the $L^2_{\mu}$-norm sense. Making use of the fact that for the low-dimensional structure $\mu$ the generalized Poincar{\'e} inequality is satisfied combining it with $\int_{\Omega}u(t)d\mu=\int_{\Omega}u^*d\mu=0$ we conclude that $u \xrightarrow{t\to \infty} u^*$ in $\mathring{H^1_{\mu}}.$ 
\qed
\end{dow}
\end{tw}

\begin{kom}
It seems possible to adapt the proof of Theorem 1 of \cite{goldie} to a wider class of differential equations on low-dimensional structures like, for instance, the p-Laplace equation with the force term dependent on a solution. Nevertheless, this kind of differential equations in the framework of low-dimensional structures has not been studied yet, and we leave it for future studies.
\end{kom}

\end{document}